\documentclass[AMA,STIX1COL,doublespace]{WileyNJD-v2}

\articletype{Research Article}%

\usepackage{amsmath}

\usepackage{makeidx}         
\usepackage{graphicx}        
\usepackage{multicol}        
\usepackage[bottom]{footmisc}
\usepackage{subfigure}
\usepackage{booktabs}

\usepackage{epsfig} 
\usepackage{epstopdf}
\usepackage{color}
\usepackage{colordvi}

\newcommand{\dualp}[2]         { \big \langle {#1},{#2} \big \rangle }       
\newcommand{\norm}[2]         { \| {#1} \|_{#2} }                      
\newcommand{\bfm}[1]             { \mathbf{#1}     }             %
\newcommand{\qq}             { \bfm{q}     }             
\newcommand{\rr}             { \bfm{r}     }             
\newcommand{\Nabla}       { \boldsymbol{\nabla} }   
\newcommand{\SCZO}            { C^0(\Omega) }                    
\newcommand{\SLTO}            { L^2(\Omega) }                       
\newcommand{\SHOO}            { H^1(\Omega) }                         
\newcommand{\SHOD}            { H^1_{\Gu}(\Omega) }                         
\newcommand{\SHOP}            { H^1(\Ph) }                              
\newcommand{\SHdivP}            { H(\textcolor{black}{\text{\bf{div}}},\Ph) }            
\newcommand{\VVK}            { {V(K_m)} }              
\newcommand{\VOM}            { {V(\Omega)} }              
\newcommand{\VV}            { {V(\Ph)} }        
\newcommand{\UUU}            { U(\Omega) }           
\newcommand{\UUUD}            { U^D(\Omega) }           
\newcommand{\UUUh}            { U^h(\Omega) }           
\newcommand{\VVh}            { V^*(\Ph) }          
\newcommand{\SHOK}            { H^1(K_m) }                        
\newcommand{\SLTK}            { L^2(K_m) }                    
\newcommand{\SHdivO}            { H(\text{\bf{div}},\Omega) }             
\newcommand{\SHdivK}            { H(\text{\bf{div}},K_m) }            
\newcommand{\SHHdK}            { H^{1/2}(\dKm) }         
\newcommand{\SHmHdK}            { H^{-1/2}(\dKm)  }               
\newcommand{\UUUhat}            { \hat{U}(\Gamma_h) }    

\newcommand{\xx}             { \bfm{x}     }             
\newcommand{\yy}             { \bfm{y}     }             
\newcommand{\err}             { {(\bfm{e},\bfm{E})}   }       
\newcommand{\errH}             { {(\bfm{e}^h,\bfm{E}^h)}   }       

\newcommand{\Ph}            { \mathcal{P}_h }
\newcommand{\Kep}           { K_m \in \Ph}
\newcommand{\dKm}           { \partial K_m }

\newcommand{\Gtt}              { {\Gamma_{t}} }
\newcommand{\Gu}              { {\Gamma_{\uu}} }

\newcommand{\dx}              { \; {\rm d} \bfm{x}   }                
\newcommand{\dss}             { \, {\rm d} s   }                          
\newcommand{\summa}[2]        { \overset{#2}{\underset{#1}{\sum}} } 
\newcommand{\supp}[1]         { \underset{#1}{\sup} \, }        

\newcommand{\wwwh}              { \bfm{w^*} }                      
\newcommand{\vvi}              { \boldsymbol{\tilde{\phi}}^i }                 
\newcommand{\wwwi}              { \boldsymbol{\tilde{\varphi}}^i }          

\newcommand{\SHR}            { H^r(\Omega) }                       
\newcommand{\SHS}            { H^s(\Omega) }                       

\newcommand{\SLinfty}                {{ L^\infty(\Omega) }  }                            
\newcommand{\SHMHGt}            { {H^{-1/2}(\Gtt) }          }                    
\newcommand{\uu}              { \bfm{u} }                      
\newcommand{\tr}               { \bfm{t}     }             
\newcommand{\ff}               { \bfm{f}     }             
\newcommand{\E}               { \bfm{E}     }             

\newcommand{\vv}              { \bfm{v} }                      
\newcommand{\zzz}              { \bfm{z} }                      
\newcommand{\www}             { \bfm{w} }                      
\newcommand{\nn}              { \bfm{n} }                      
\newcommand{\pp}              { \bfm{p} }       
\newcommand{\sig}             { \boldsymbol{\sigma}}           
\newcommand{\eps}             { \boldsymbol{\varepsilon} }     
\newcommand{\isdef}           { \overset{\text{def}}{=} } 
\newcommand{\ds}              { \displaystyle }   

\theoremstyle{plain}
\newtheorem{thm}{Theorem}[section]
\newtheorem{lem}{Lemma}[section]
\newtheorem{rem}{Remark}[section]

\newtheorem{prp}{Proposition}[section]

\begin{document}

\title{A Stable Mixed FE Method for Nearly Incompressible Linear Elastostatics}

%
%

\author[1]{Eirik Valseth*}

\author[2]{Albert Romkes}

\author[2]{Austin R. Kaul}

\author[1]{Clint Dawson}

\authormark{EIRIK VALSETH \textsc{et al}}

\address[1]{\orgdiv{Oden Institute for Computational Engineering and Sciences}, \orgname{The University of Texas at Austin}, \orgaddress{\state{Texas}, \country{USA}}}

\address[2]{\orgdiv{Department of Mechanical Engineering}, \orgname{South Dakota School of Mines \& Technology}, \orgaddress{\state{South Dakota}, \country{USA}}}

\corres{*Eirik Valseth, Oden Institute for Computational Engineering and Sciences, The University of Texas at Austin, Austin, TX 78712, USA. \email{eirik@utexas.edu}}





\abstract[Summary]{
We present a new, stable,  mixed finite element (FE) method for linear elastostatics of nearly incompressible solids. The method is the automatic variationally stable FE (AVS-FE) method of Calo, Romkes and Valseth,
in which we consider a Petrov-Galerkin weak formulation where the stress and displacement variables are in the space $H(\bf{div})$ and $H^1$, respectively.
This allows  us to employ a fully conforming FE discretization for any elastic solid using classical FE subspaces of $H(\bf{div})$ and $H^1$. Hence, the resulting FE approximation yields both continuous stresses and displacements. 

To ensure stability of the method, we employ the philosophy of the discontinuous Petrov-Galerkin (DPG) 
method of 
Demkowicz and Gopalakrishnan  and use optimal test spaces. 
Thus, the resulting FE discretization is stable even as the Poisson's ratio $\nu \rightarrow 0.5$, and the system of linear algebraic equations is symmetric and positive definite.
Our method also comes with a built-in \emph{a posteriori} error estimator as well as indicators 
which are used to drive mesh adaptive refinements. 
We present several numerical verifications of our method including comparisons to existing FE technologies.
}
\maketitle

\section{Introduction}
\label{sec:introduction}

Linear elastostatics is arguably the most successful area of application of the classical (Bubnov-Galerkin) FE
method. For homogeneous isotropic engineering materials, such as steel and aluminum, the Bubnov-Galerkin 
method is stable, satisfies a best approximation property in terms of elastic strain energy and is 
computationally efficient.  However, for commonly employed modern engineering materials such as rubbers 
and soft plastics, i.e., nearly incompressible materials, the Bubnov-Galerkin method suffers from
\emph{locking} and loss of discrete stability 
(see, e.g.,~\cite{babuvska1992lockingi,babuvska1992locking,oden2012introduction}). 
In~\cite{phillips2009overcoming}, Phillips and Wheeler investigate this phenomenon in great detail and
highlight that error estimates for the classical FE method depend on the factor $1/(1-2\nu)$, which
clearly tends to infinity as $\nu \rightarrow 0.5$.


Mixed FE methods~\cite{BrezziMixed} provide functional settings in which certain FE discretizations 
are stable for mixed forms of the elastostatics problem as well as for the Stokes 
equations which and can be shown to be equivalent to the equations of linear 
elastostatics when $\nu = 0.5$.
Other mixed FE methods based on the consideration of the underlying equations of elastostatics using 
the compliance tensor do not suffer from this loss of stability, but often require additional  
constraints to ensure symmetric stresses~\cite{falk2008finite}, called Hellinger-Reissner formulations, 
and are typically more computationally costly than the primal Bubnov-Galerkin formulation.
In~\cite{arnold2002mixed}, Arnold and Winther present an element for the Hellinger-Reissner formulation
which was one of the first stable elements utilizing polynomial bases for both stress and displacement.
The difficulty in building conforming approximation spaces for the stress has also resulted in several
nonconforming mixed methods, see e.g.,~\cite{arnold2003nonconforming,gopalakrishnan2012second} where 
the FE approximation of the stress does not reside in the space dictated by the weak formulation. 
The task of establishing approximation spaces for these mixed FE methods is certainly not trivial 
and is an active area of investigation and 
recent publications include~\cite{quinelato2019full,ambartsumyan2020multipoint}. 

Stabilized FE methods that adjust the functionals of the weak formulation can be used to ensure 
discrete stability~\cite{chiumenti2002stabilized}. This type of stabilization is performed for
 both the mixed and classical FE methods, see the work of Nakshatrala 
\emph{et al.}~\cite{nakshatrala2008finite} as well as Masud \emph{et al.}~\cite{masud2005stabilized}.
However, stabilized methods generally require 
arduous analyses to establish a proper choice of penalization/stabilization parameters. 
Reduced integration methods are also commonly used when approaching the incompressible
limit~\cite{malkus1978mixed}.
The discontinuous  Galerkin (DG) method also remains a popular choice for nearly incompressible 
elastostatics~\cite{hansbo2002discontinuous,phillips2009overcoming,liu2009three}. \textcolor{black}{In general}, these 
achieve stability by adjusting the inter-element jump or average terms by weights in a manner similar to 
the stabilized FE methods.   

Stable FE methods such as the least squares FE methods (LSFEMs) (see, e.g., text by Bochev and 
Gunzberger~\cite{bochevLeastSquares}) or the discontinuous Petrov-Galerkin (DPG) Method of 
Demkowicz and Gopalakrishnan~\cite{Demkowicz5} can be employed to resolve the stability issue. 
The LSFEM has been applied to linear elastostatics in~\cite{cai2003first}
and in~\cite{cai2005adaptive}, a weighted first-order system least squares is applied successfully 
to nearly incompressible materials.
Gopalakrishnan and Qiu  provide an analysis of the well-posedness of the DPG method applied to linear 
elastostatics  in~\cite{gopalakrishnan2014analysis}.
In~\cite{bramwell2012locking}, Bramwell \emph{et al.} consider two distinct DPG methods for the 
nearly incompressible elastostatics problem that are locking free and present numerical verifications
highlighting capabilities as $\nu \rightarrow 0.5$. 
 The DPG has also been successfully 
applied to this problem in several works, including the fully incompressible case in~\cite{fuentes2017coupled} employing the compliance tensor to avoid locking in that case.
In~\cite{fuentes2017coupled,keith2016dpg}, the DPG method is applied to the problem  
of linear elastostatics and several variational formulations are considered 
including for the case of nearly incompressible materials. 
In particular, in~\cite{fuentes2017coupled}, the idea of coupling multiple weak formulations 
throughout the computational domain is explored in great detail.

In the classical FE method, the approximations of displacements of the equivalent weak form of the 
underlying partial differential equation (PDE) of static equilibrium are sought in $C^0$ continuous
polynomial spaces and stress approximations are established by computing gradients of the displacements,
 i.e., the stresses are piecewise discontinuous.
On the other hand, mixed FE methods for the linear elastostatics problem consider an equivalent 
first-order system of the underlying PDE.
This first-order system description can lead to weak forms which allow stresses to be in $\SHdivO$
and displacements that are in $\SLTO$ (see Section 2.4 of~\cite{keith2016dpg} for a thorough 
discussion on other options). Hence, in the FE approximations the displacements must be sought in 
piecewise discontinuous polynomial function spaces.
The theory of distributions ensures that optimally convergent FE solutions can be established
for both classical and mixed FE methods, as well as for their properly stabilized counterparts if
$\nu$ is close to $0.5$.
However, the resulting numerical approximations are not physical, as we know that both the displacement and certain components of the stress fields are continuous.
We know of \textcolor{black}{three options} to establish both continuous displacements and stresses. 
\textcolor{black}{Firstly}, the isogeometric 
FE methods of Hughes \emph{et al.}~\cite{hughes2005isogeometric} which uses higher order bases for the 
discrete FE approximation, i.e., $C^k$ continuity. \textcolor{black}{Secondly},  the $k-$version FE method of Surana 
\emph{et al.}~\cite{Surana2010} which employs higher order bases as well as a least squares approach.
The  popularity of the isogeometric FE method has grown significantly over the last decade, but the 
stability issue of nearly incompressible materials 
still persist. In~\cite{taylor2011isogeometric}, Taylor introduces a mixed version of the isogeometric 
FE methods for incompressible solids where discontinuous stress approximations are sought.  
\textcolor{black}{Lastly, the use of post processing techniques where a discontinuous solution component is projected into a continuous discrete space, e.g., by using Oswald operators, 
see~\cite{oswald1997intergrid,ainsworth2005robust} for details and further references.}


The automatic variationally stable finite element (AVS-FE) method  
introduced by Calo, Romkes and Valseth in~\cite{CaloRomkesValseth2018} provides a framework, 
much like the DPG of Demkowicz and Gopalakrishnan~\cite{Demkowicz5}, to establish stable FE approximations for any PDE. However, the AVS-FE differs in its choice of trial spaces while 
employing the DPG concept of optimal discontinuous test functions. 
In addition to the approximation of the trial variables, the AVS-FE also comes with a "built-in" error estimator and error indicators that can be employed to drive mesh adaptive strategies.
The stability property of the AVS-FE allows us to derive Petrov-Galerkin weak formulations that are
posed with trial functions that are in classical Hilbert spaces, e.g., $H(\bf{div})$ and $H^1$.
Hence, the corresponding FE approximations are to be sought in classical \emph{continuous} FE approximation spaces yielding \emph{continuous} FE approximations for all trial variables. 
The LSFEMs presented in~\cite{cai2003first,cai2005adaptive} also pose weak formulations in Hilbert 
spaces as the AVS-FE but considers alternative formulations for the elasticity problem and considers 
nonconforming approximations for the displacement.

In this paper, we build upon the preliminary investigation of Valseth in~\cite{eirik2019thesis} 
for the AVS-FE method applied to linear elastostatics of nearly incompressible media.
We introduce our model problem and notations in Section~\ref{sec:model_and_notation},
The weak formulation and its corresponding FE discretzation are presented in Section~\ref{sec:avs-fe}
in conjunction with a brief review of the AVS-FE methodology.
In Section~\ref{sec:error_estimates}, we present optimal \emph{a priori} error estimates.
Several numerical verifications are presented in Section~\ref{sec:experiments} highlighting 
the stability of our method as $\nu \rightarrow 0.5$, including an asymptotic convergence 
study with comparisons to existing FE methods. 
We draw conclusions and discuss future works in Section~\ref{sec:conclusions}.

\section{The AVS-FE Method}  
\label{sec:avs-fe}
The AVS-FE method~\cite{CaloRomkesValseth2018} provides  
a functional setting to analyze linear boundary value problems (BVPs) in which the
underlying differential operator is non self-adjoint or leads to unstable FE discretizations. 
In this section we introduce our model problem, and briefly review the AVS-FE method. 
A thorough introduction can be found in ~\cite{CaloRomkesValseth2018}.

\subsection{Model Problem: Linear Elastostatics of Nearly Incompressible Solids}
\label{sec:model_and_notation}
\begin{figure}[t]
\centering
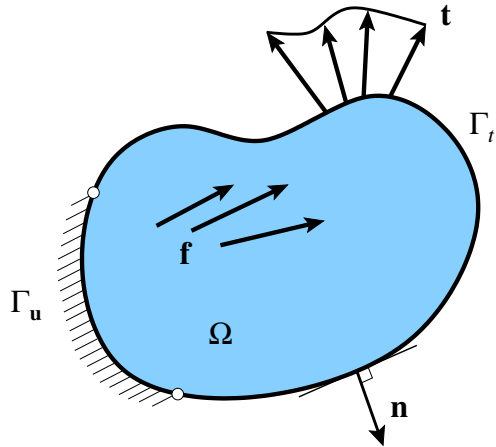
\caption{The model problem.}
\label{fig:model}
\end{figure}

Let $\Omega\subset \mathbb{R}^{\textcolor{black}{n}}$, \textcolor{black}{$n=1,2,3$ (we consider the two dimensional case here for simplicity)} be a bounded open domain,  containing a linearly elastic, nearly incompressible, and possible heterogeneous solid. The boundary  $\partial \Omega$ is partitioned into two open and disjoint segments $\Gtt$ and $\Gu$, such that $\partial \Omega = \overline{ \Gtt \cup \Gu}$. As depicted in Figure~\ref{fig:model},  the body is in static equilibrium under the action of external  body loads $\ff\in[\SLTO]^2$ in $\Omega$, surface tractions  $\tr\in\SHMHGt$ on $\Gtt$, as well \textcolor{black}{as fixed zero} displacements on $\Gu$. Since the solid is assumed to be linearly elastic, its constitutive behavior is governed by  Generalized Hooke's Law, i.e.:
\begin{equation}
\label{eq:Hooke}
  \sig = \E \, \eps,
\end{equation}
where $\sig$ denotes the (2D) Cauchy stress tensor, $\eps$ the (2D) Green strain tensor, and  $\E$ the fourth order (Riemann) elasticity tensor, with elliptic and symmetric Riemann coefficients $E_{ijkl}\in\SLinfty$. In this work, we limit our focus to problems in which the deformations in the material remain small and therefore the kinematic relation between the strain tensor $\eps$ and displacement field $\uu$ is linear and governed by:
%
\begin{equation}
\label{eq:kinematic_relation}
  \eps = \frac{1}{2} \left[ \Nabla \uu + (\Nabla \uu)^T \right].
\end{equation}
With these notations and relations in force, the equilibrium state of the solid is represented by the following BVP, governing the displacement field $\uu$:
\begin{equation}
\label{eq:model_pde}
  \boxed{
  \begin{array}{l}
      \text{Find} \;  \uu\in\textcolor{black}{[\SHOO]^2} \;  \text{such that:}   \\
   \\[-0.1in]
     \hspace{0.5in} \begin{array}{rcll}
        - \Nabla \cdot  \textcolor{black}{\sig}  &  = & \ff, &
\text{in} \; \Omega,
        \\[0.1in]
        \textcolor{black}{\sig} \nn &  = & \tr,        & \text{on} \; \Gtt,
       \\
         \uu              &  = & \bfm{0},     & \text{on} \; \Gu,
      \end{array}
  \end{array}
  }
\end{equation}
where $\nn$  denotes the outward unit normal vector to $\partial \Omega$.  In this paper, we consider the specific scenario in which the  solid is comprised of one or more constituents with nearly incompressible material properties. Hence, the Riemann coefficients $E_{ijkl}$ can involve values of the Poisson Ratio  $\nu$ that are very close to, but still less than, $0.5$.

In the following, we shall use the following notations:
\begin{itemize}

\item inner products between vector
valued functions are denoted with the single dot symbol $\cdot$, 
and inner products between tensor valued functions are denoted by the colon or double dot symbol $\colon$. 

\item $h_m$ is the diameter  of element $K_m$.

\item in weak formulations, we present edge integrals using trace the operators: $i)$ $\gamma^m_0: \SHOK: \longrightarrow \SHHdK$ as the local zeroth order trace operator and $ii)$ $\gamma^m_{\nn_m}:\SHdivK \longrightarrow \SHmHdK$ denote the local normal trace operators where $\nn_m$ is the outward unit normal vector to the element boundary $\dKm$ (e.g., see~\cite{Girault1986}).
\item \textcolor{black}{vector and tensor valued test functions are denoted using  $\vv$ and $\www$, respectively. Restrictions of these to an element $K_m$ are denoted by employing the subscript $m$.}

\end{itemize}

\subsection{Weak Formulation}
\label{sec:weak_form}

AVS-FE weak formulations are established using techniques that are similar to DG and DPG methods by 
considering element-wise weak formulations that are subsequently summed throughout the FE mesh 
partition to yield global weak formulations. We mention only key points here and omit the full derivation here for brevity but refer 
to~\cite{CaloRomkesValseth2018} for detailed derivations. 

To establish AVS-FE weak formulations, we first require a partition
$\Ph$ of $\Omega$ into convex elements $K_m$, such that:
\begin{equation}
\notag
\label{eq:domain}
  \Omega = \text{int} ( \bigcup_{\Kep} \overline{K_m} ), \quad K_m \cap K_n \textcolor{black}{=0}, \quad m \ne n.
\end{equation}
The partition $\Ph$ is such that any discontinuities in $E_{ijkl}$ are
restricted to the boundaries of each element $\dKm$.
The BVP~\eqref{eq:model_pde} is recast as a first-order system by using the stress tensor 
from the constitutive law~\eqref{eq:Hooke}, i.e.,
\begin{equation}
\label{eq:model_pde_first_order}
  \boxed{
  \begin{array}{l}
      \text{Find} \;  (\uu,\sig)\in\textcolor{black}{[\SHOO]^2}\times\SHdivO \;  \text{such that:}   \\
   \\[-0.1in]
     \hspace{0.5in} \begin{array}{rcll}
      \sig - \E \, \eps &  = & \bfm{0} &
\text{in} \; \Omega,
        \\
        - \Nabla \cdot \sig &  = & \ff, &
\text{in} \; \Omega,
        \\[0.1in]
        \sig \, \nn  &  = & \tr,        & \text{on} \; \Gtt,
       \\
         \uu              &  = & \bfm{0},     & \text{on} \; \Gu,
      \end{array}
  \end{array}
  }
\end{equation}
where $\SHdivO$ is the Hilbert space of tensor-valued functions which divergence is weakly continuous
and $\eps = \eps(\uu)$ denotes the gradient operator in~\eqref{eq:kinematic_relation}. 
Note that this first-order system, or mixed form, BVP is standard for mixed FE methods for
linear eleastostatics.

Next, the first-order system is multiplied by test functions $(\vv,\www)\in\SLTK^6$ and enforced 
weakly on each individual element $\Kep$. We then apply integration by parts locally on each element $K_m$ to the term involving the 
divergence of the stress field $\Nabla \cdot \sig$ to enable weak applications of the Neumann 
boundary condition (BC). A subsequent summation of all elements in $\Ph$
and strong enforcement of \textcolor{black}{all} BCs leads us to the AVS-FE weak formulation: 
\begin{equation} \label{eq:weak_form}
\boxed{
\begin{array}{ll}
\text{Find } (\uu,\sig) \in \UUU \text{ such that:}
\\[0.05in]
 \qquad   B((\uu,\sig),(\vv,\www)) = F(\vv), \quad \forall (\vv,\www)\in \VV, 
 \end{array}}
\end{equation}
where the test space $\VV$ is broken, the  bilinear form, $B:\UUU\times\VV\longrightarrow \mathbb{R}$, and linear functional, $F:\VV\longrightarrow \mathbb{R}$ are defined:
\begin{equation} \label{eq:B_and_F_LE}
\begin{array}{c}
B((\uu,\sig),(\vv,\www)) \isdef
\ds \summa{\Kep}{}\biggl\{ \int_{K_m} \biggl[  \, \left(  \sig - \E \eps(\uu)\right) \textcolor{black}{:} \www_m \, 
+  \, \sig : \Nabla \vv_m  \biggr] \dx \biggr.
\oint_{\dKm}  \gamma^m_\nn(\sig) \, \gamma^m_0(\vv_m) \, \dss \biggr\},
 \\[0.15in]
 F(\vv) \isdef   \ds \summa{\Kep}{} \int_{K_m} \ff \cdot \vv_m \dx , 
 \end{array}
\end{equation}
%
and the function spaces are defined:
\begin{equation}
\label{eq:function_spaces_LE}
\begin{array}{c}
\UUU \isdef \biggl\{ (\uu,\sig)\in \textcolor{black}{[\SHOO]^2}\times\SHdivO: \; \gamma_0^m(\uu)_{|\dKm\cap\Gu} =\bfm{0}, \textcolor{black}{\gamma_0^m(\sig)_{|\dKm\cap\Gtt} =\tr } \;  \forall\Kep\biggr\},
\\[0.15in]
\VV \isdef \biggl\{ (\vv,\www)\in [\SHOP]^2\times[\SLTO]^4: \,  \gamma_0^m(\vv_m)_{|\dKm\cap\Gu} = \bfm{0}, \; \forall\Kep\biggr\},
\end{array}
\end{equation}
with norms $\norm{\cdot}{\UUU}:  \UUU \!\! \longrightarrow \!\! [0,\infty)$ and $\norm{\cdot}{\VV}: \VV\! \! \longrightarrow\! \! [0,\infty)$ defined as:
\begin{equation}
\label{eq:norms}
\begin{array}{l}
\ds \norm{(\uu,\sig)}{\UUU} \isdef \sqrt{\int_{\Omega} \biggl[ \Nabla \uu : \Nabla \uu + \uu \cdot \uu   + (\Nabla \cdot \sig)^2+\sig \cdot \sig\biggr] \dx },
\\[0.2in]
\ds   \norm{(\vv,\www)}{\VV} \isdef \sqrt{\summa{\Kep}{}\int_{K_m} \biggl[  h_m^2 \Nabla \vv_m : \Nabla \vv_m + \vv_m \cdot \vv_m   + \www_m : \www_m\biggr] \dx }.
 \end{array}
\end{equation}
The norm $\norm{\cdot}{\VV}$ is equivalent to the \textcolor{black}{$L^2$} norm:
\begin{equation}
\label{eq:norm2}
\begin{array}{l}
\hspace{-0.5cm} \ds   \norm{(\vv,\www)}{\SLTO} \textcolor{black}{=} \sqrt{ \int_{\textcolor{black}{\Omega}} \biggl[ 
  \textcolor{black}{\vv}  \cdot \textcolor{black}{\vv} + \www : \www  \biggr] \dx }.
 \end{array}
\end{equation}
Note that the edge integrals in~\eqref{eq:B_and_F_LE} are to be interpreted as duality pairings 
in $\SHHdK\times\SHmHdK$, but we employ notation that is engineering convention here using the 
integral representation. Most importantly, since $(\uu,\sig)\in\textcolor{black}{[\SHOO]^2}\times\SHdivO$, these integrals 
are well defined and our trial space is continuous. 
As the trial and test spaces are of different regularity we are in a Petrov-Galerkin setting, 
particularly a DPG setting, since the test space is broken. However, since the trial space is
$\textcolor{black}{[\SHOO]^2}\times\SHdivO$ our functional setting differs from that of DPG methods in which the regularity of the trial space is reduced by introducing variables on the edge of each element. 
\begin{rem}
\textcolor{black}{The norm $\norm{\cdot}{\VV}$ is used in the discrete computational setting due to exhaustive numerical experimentation. In particular, its use is justified based on $i)$ engineering and  $ii)$ computational intuition. $i)$ The scaling ensures consistency of units, e.g., if $\vv$ is of unit length then $\Nabla \vv$ is of unit $\frac{1}{\text{length}}$. Thus, all entries are of the same unit. $ii)$ The scaling also ensures that all terms in the norm are of the same magnitude in the discrete setting due to the exact same reasoning since gradient terms scale as $h^{-1}$. Finally we note that the equivalence between $\norm{\cdot}{\VV}$ and $\norm{\cdot}{\SLTO}$ is based on mesh dependent constants. 
}
\end{rem}

In the spirit of the DPG method, we now introduce an equivalent norm on the trial space, the energy norm $\norm{\cdot}{\text{B}}: \UUU\longrightarrow [0,\infty)$:
\begin{equation}
\label{eq:energy_norm}
\norm{(\uu,\sig)}{\text{B}} \isdef \supp{(\vv,\www)\in \VV\setminus \{(\mathbf{0},\mathbf{0})\}} 
     \frac{|B(\uu,\sig),(\vv,\www))|}{\quad\norm{(\vv,\www)}{\VV}},
\end{equation}
and a Riesz representation problem for $(\pp,\rr)$, the optimal test functions:
\begin{equation}
\label{eq:riesz_problem}
\begin{array}{rcll}
\ds \left(\, (\pp,\rr),(\vv,\www) \, \right)_\VV &  \! \! =  \! & B(\,(\uu,\sig),(\vv,\www) \, ),& \, \forall (\textcolor{black}{\vv} ,\www)\in\VV. 
\end{array}
\end{equation}
The Riesz representation problem is well posed with unique solutions due to the inner product in the left hand side (LHS) and guarantees the stability of DPG methods. 
The Riesz representation problem also leads to the following norm equivalence:
\begin{equation}
\label{eq:norm_equivalence}
\norm{(\uu,\sig)}{\text{B}} = \norm{(\pp,\rr)}{\VV},
\end{equation}
which will be employed extensively in the following. 
For details on optimal test functions and proof of the norm equivalence, we refer 
to~\cite{Demkowicz2,Demkowicz5}.
Due to the energy norm,
the bilinear form~\eqref{eq:B_and_F_LE} has continuity and \emph{inf-sup} constants equal to one and 
the load functional is also continuous which can be shown using classical techniques. Hence, we have 
established a well posed weak formulation of the linear elastostatics BVP using continuous trial spaces 
for both displacement and stress fields, i.e., $\textcolor{black}{[\SHOO]^2}\times\SHdivO$, in terms of the energy 
norm~\eqref{eq:energy_norm}.

Well-posedness in terms of the energy norm is essentially an assumption of DPG methods as we define 
a norm that ensure \emph{inf-sup} and continuity conditions of the bilinear form. For completeness, 
we also provide the following lemma of well-posedness in standard Sobolev norms by first stating 
two important results. \textcolor{black}{For the sake of simplicity, we consider the case in which homogeneous Dirichlet BCs are applied on the full boundary $\partial \Omega = \Gu$}:
\begin{prp} \label{prp:dual_wellp}
Let $(\uu^{\rm{D}}, \sig^{\rm{D}}) \in [\SHOD]^2\times[\SLTO]^4$ be the solution of the dual mixed
formulation:
\begin{equation} \label{eq:B_F_dual}
\begin{array}{ll}
\text{Find } (\uu^{\text{D}}, \sig^{\text{D}}) \in [\SHOD]^2\times[\SLTO]^4 \text{ such that:}
\\[0.05in]
\ds \quad  \underbrace{ \int_{\Omega} \biggl[ \left(  \sig^{\text{D}} - \E \eps(\uu^{\text{D}})\right) \textcolor{black}{:} \www
  +   \sig^{\text{D}} : \Nabla \vv  \biggr] \dx }_{b^{\text{D}}(\uu^{\text{D}}, \sig^{\text{D}}),(\vv,\www))} = \int_{\Omega}  \ff \cdot \vv \dx , \quad \forall (\vv,\www)\in [\SHOD]^2\times[\SLTO]^4,
\end{array}
\end{equation}
which is well posed. Hence, the bilinear form satisfies the 
\emph{inf-sup} condition:
\begin{equation} \label{eq:inf-supp}
\begin{array}{l}
\exists \gamma >0 : \supp{(\vv,\www) \in \VOM } \ds 
     \frac{|b^{\text{D}}(\uu^{\text{D}}, \sig^{\text{D}}),(\vv,\www))|}{\norm{(\vv,\www)}{\VOM}} \ge \gamma \, \norm{(\uu^{\text{D}}, \sig^{\text{D}})}{\UUUD},
 \end{array}
\end{equation}
where $\UUUD=\VOM=[\SHOD]^2\times[\SLTO]^4$ and $\SHOD$ is the space of $H^1$ functions that satisfy homogeneous Dirichlet conditions on $\Gu = \partial \Omega$.
\end{prp}
\emph{Proof}: see Theorem 2.1 in~\cite{keith2016dpg}.
\newline \noindent ~\qed

 We write the bilinear
form~\eqref{eq:B_and_F_LE} as: 
\begin{equation} \label{eq:AVS-weak_mod}
\begin{array}{l}
B((\uu,\sig),(\vv,\www)) = b^{\text{D}}((\uu,\sig),(\vv,\www)) + \dualp{\gamma^m_\nn(\sig)  }{\gamma^m_0(\vv)}_{\Gamma_h},
 \end{array}
\end{equation}
where $\dualp{\gamma^m_\nn(\sig)  }{\gamma^m_0(\vv)}_{\Gamma_h} \isdef \summa{\Kep}{}
\oint_{\dKm} \{ \gamma^m_\nn(\sig) \,, \gamma^m_0(\vv_m)  \,  \} \dss$. 
\begin{prp} 
\label{prp:zero_skel_trac}
Let $\sig \in \SHdivP$ and $\textcolor{black}{\vv} \in \textcolor{black}{[\SHOP]^2}$. Then:
\begin{equation} \label{eq:inf_sup_skel}
\begin{array}{l}
\exists \gamma^S >0 : \supp{\textcolor{black}{\vv_m} \in \textcolor{black}{[\SHOP]^2} } \ds 
     \frac{|\dualp{\gamma^m_\nn(\sig)  }{\gamma^m_0(\vv_m)}_{\Gamma_h}|}{\norm{\textcolor{black}{\vv}}{\SHOP}} \ge \gamma^S \, \norm{\sig}{\UUUhat},
 \end{array}
\end{equation}
where $\SHdivP$ denotes the broken $H(\bf{div})$ space and $\norm{\sig}{\UUUhat}$ is the minimum 
energy extension norm:
\begin{equation} \label{eq:min_en_ext}
\begin{array}{l}
\norm{\sig}{\UUUhat}  \isdef 
\ds \supp{\textcolor{black}{\vv} \in \textcolor{black}{[\SHOP]^2} } \ds 
     \frac{|\dualp{\gamma^m_\nn(\sig)  }{\gamma^m_0(\vv_m)}_{\Gamma_h}|}{\norm{\textcolor{black}{\vv}}{\SHOP}} = \,\rm{inf} \, \norm{\sig}{\SHdivO}.
 \end{array}
\end{equation}
Additionally, if $\sig \in \SHdivO$:
\begin{equation} \label{eq:zero_traces}
\begin{array}{l}
\dualp{\gamma^m_\nn(\sig)  }{\gamma^m_0(\vv)}_{\Gamma_h} = 0, \quad \textcolor{black}{\forall  \vv  \in H^{1}(\Omega)}.
 \end{array}
\end{equation}
\end{prp} 
\emph{Proof}: see Theorem 2.3 in~\cite{carstensen2016breaking}.
\newline \noindent ~\qed
\begin{lem}
Let $(\uu, \sig) \in \UUU$  and $(\vv, \www) \in \VV$. Then, the AVS-FE weak
formulation~\eqref{eq:weak_form} satisfies all conditions of the Babu{\v{s}}ka Lax-Milgram Theorem~\cite{babuvska197finite} and is well posed.
\end{lem}
\emph{Proof}: The load functional and bilinear form~\eqref{eq:B_and_F_LE} are continuous due to
the Cauchy-Schwarz inequality. 
The following \emph{inf-sup} condition:
\begin{equation} \label{eq:inf-sup}
\begin{array}{l}
\exists C >0 : \supp{(\vv,\www) \in \VV } \ds 
     \frac{|B(\uu, \sig),(\vv,\www))|}{\norm{(\vv,\www)}{\textcolor{black}{W(\Ph)}}} \ge \gamma \, \norm{(\uu, \sig)}{\UUU},
 \end{array}
\end{equation}
is satisfied due to Theorem 3.3 in~\cite{carstensen2016breaking}. \textcolor{black}{This theorem holds if
the following conditions hold (see Assumptions 3.1 and 3.2 in~\cite{carstensen2016breaking}):
$i)$ the bilinear form $b^{\text{D}}(\cdot,\cdot)$ 
satisfies the \emph{inf-sup} condition and has a trivial kernel, the form $\dualp{\gamma^m_\nn(\sig)  }{\gamma^m_0(\vv)}_{\Gamma_h}$  satisfies $ii)$ an \emph{inf-sup} condition and  $iii)$ a kernel preserving property. Due to Proposition~\ref{prp:dual_wellp}, we can conclude that the bilinear form  $b^{\text{D}}(\cdot,\cdot)$ satisfies $i)$. Second, $\dualp{\gamma^m_\nn(\sig)  }{\gamma^m_0(\vv)}_{\Gamma_h}$  satisfies the inf-sup condition in~\eqref{eq:inf_sup_skel}, and the kernel preserving property is satisfied by noting that $\dualp{\gamma^m_\nn(\sig)  }{\gamma^m_0(\vv)}_{\Gamma_h}$ vanishes if evaluated using test functions from the test space of $b^{\text{D}}(\cdot,\cdot)$.
\textcolor{black}{The norm in the denominator is defined: $\norm{(\vv,\www)}{\textcolor{black}{W(\Ph)}}^2 = \norm{\vv}{\SHOP}^2+\norm{\www}{\SLTO}^2 $. } }
%
\newline \noindent ~\qed

\begin{rem} \label{rem:alternative_forms}
The bilinear and linear forms in~\eqref{eq:B_and_F_LE} are not unique choices for the AVS-FE method. 
We have chosen these particular forms as they allow us to keep the weak formulation close to classical 
mixed FE methods for linear elastostatics and enforce Dirichlet BCs strongly and Neumann BCs weakly. 
Other forms can be derived in which the 
trial space is continuous, and the test space is discontinuous, due to the flexibility
of the Petrov-Galerkin method.  In~\cite{keith2016dpg} Keith~\emph{et al.} consider several possible weak formulations for the elastostatics problem and the DPG method and perform a rigorous analysis showing their well posedness. 
\end{rem}
%
%

%
%
%
%
%

%

\subsection{AVS-FE Discretization}
\label{sec:discretization}

To establish FE discretizations of the weak formulation~\eqref{eq:weak_form}, the AVS-FE takes the 
approach of classical FE methods and seeks continuous polynomial approximations that are in conforming
subspaces of the $\textcolor{black}{[\SHOO]^2}\times\SHdivO$ trial spaces. Hence, for the displacement field we use classical $\SCZO$ continuous Lagrange polynomials. Generally, in mixed
FE methods this choice leads to unstable and inconsistent FE discretizations and is avoided and the stress field
is sought in a Raviart-Thomas (RT) or Brezzi-Douglas-Marini (BDM) space~\cite{BrezziMixed}. 
Due to the stability properties of the AVS-FE method, it is often convenient  
to employ the same $\SCZO$ polynomials for the stress variable. As reported in~\cite{valseth2020goal}, 
for convex domains and smooth solutions, the $\SCZO$ are superior. 
In Section~\ref{sec:experiments} we present numerical verifications comparing the approximations from 
these spaces. Furthermore, we present a verification where we again compare the classical RT spaces with the $\SCZO$ polynomials for a physical application in which the stress 
field is such that it is discontinuous in the tangential direction.

Hence, we seek numerical approximations $(\uu^h,\sig^h)$ of $(\uu,\sig)$ of the
weak formulation~\eqref{eq:weak_form} and  represent the approximations as linear 
combinations of the trial bases $(\boldsymbol{\phi}^i(\xx),\boldsymbol{\varphi}^j(\xx))\in\UUUh$ (e.g., $\textcolor{black}{ \mathcal{P}^p (\Omega) \times RT_p(\Omega)}$ ) and their corresponding degrees of freedom:
\begin{equation} \label{eq:FE_sol}
\uu^h(\xx) = \summa{i=1}{N_{\uu}} \uu^h_i \, \boldsymbol{\phi}^i(\xx), 
\quad \sig^h(\xx) = \summa{j=1}{N_{\sig}} \sig^{h,j} \, \boldsymbol{\varphi}^j(\xx).
\end{equation}
Now, the test space, which is discontinuous, is to be constructed by the DPG philosophy using 
optimal test functions defined by the discrete equivalent of the Riesz representation problem~\eqref{eq:riesz_problem}. Thus, the optimal test space is spanned by functions that are solutions 
of the global weak problems, e.g., for a trial basis function $\boldsymbol{\phi}^i(\xx)$ for the displacement variable, 
its corresponding optimal test function $(\vvi,\wwwi)$ is defined by:
\begin{equation}
\label{eq:test_problem}
\begin{array}{rcll}
\ds \left(\, (\rr,\zzz),(\vvi,\wwwi) \, \right)_\VV &  \! \! =  \! & B(\,(\boldsymbol{\phi}^i,\mathbf{0}),(\rr,\zzz) \, ),& \, \forall (\rr \textcolor{black}{,\zzz)} \in \VV, \quad i=1,\dots, N.
\end{array}
\end{equation}
Inspection of~\eqref{eq:test_problem} reveals the ingenuity of the DPG philosophy since 
the test space is broken, we do not need to solve this problem globally but rather element-wise local analogues which can be solved in a completely decoupled local fashion. Hence, the AVS-FE 
approximation is governed by:
\begin{equation} \label{eq:FE_disc_weak}
\boxed{
\begin{array}{ll}
\text{Find } (\uu^h,\sig^h) \in \UUUh \text{ such that:}
\\[0.05in]
 \qquad   B((\uu^h,\sig^h),(\vv^{*},\wwwh)) = F(\vv^{*}), \quad \forall (\vv^{*},\wwwh)\in \VVh, 
 \end{array}}
\end{equation}
where the test space $\VVh$ is spanned by the approximated optimal test functions computed from 
local equivalents of~\eqref{eq:test_problem}.

The choice we have made of fully continuous trial spaces has several important consequences: $i)$ as 
the bilinear form~\eqref{eq:B_and_F_LE} is such that information is transferred from element-to-element 
by the continuity of trial functions alone, the optimal test functions in the AVS-FE have the same 
support as its trial basis functions. Hence, the resulting global stiffness matrix 
has the same sparseness as mixed FE methods.  $ii)$ the local optimal test function 
problems can be solved by using the same polynomial degree of approximation as the trial functions 
which define each problem. Thus, the cost incurred to establish the optimal test functions is kept 
as low as possible (see Remark~\ref{rem:discrete_Stab}). $iii)$ finally, the AVS-FE optimal test functions can be implemented in 
legacy FE software in which continuous polynomials are the only available basis functions by 
redefining the element stiffness matrix assembly process.
\begin{rem} \label{rem:discrete_Stab}
If the computation of the optimal test functions could be performed exactly, discrete inf-sup 
constant would be identical to the continuous one (often referred to as the ideal DPG method).
Unless the test space is $L^2$,  i.e., the LSFEM, this is not possible in practical computations 
and we consider an approximation of these functions~\cite{gopalakrishnan2014analysis}. Thus, there is a potential 
loss of discrete stability if the optimal test functions are computed without sufficient accuracy. 
Sufficient accuracy is ensured by  the existence of (local) Fortin
operators~\cite{boffi2013mixed}. The construction of such operators for the DPG method is studied in
great detail in 
~\cite{nagaraj2017construction}, and its analysis was recently further refined 
in~\cite{demkowicz2020construction}. For second order PDEs, a Fortin 
operator's existence and thus discrete stability is ensured if the local Riesz representation 
problems are solved using polynomials of order $r=p+\Delta p$, where $p$ is the degree of the trial 
space discretization and $\Delta p = d$ the space dimension. However, while this enrichment degree
ensures the existence of the required Fortin operator, numerical evidence suggest that in most cases
$\Delta p = 1$ is typically sufficient~\cite{demkowicz2020construction}. Alternative test 
spaces for the DPG method for singular perturbation problems are investigated 
in~\cite{salazar2019alternative}, even for the case of $\Delta p = 0$.

In  the AVS-FE method, numerical evidence suggest that $r=p$ is 
sufficient~\cite{CaloRomkesValseth2018,valseth2020goal} for convection-diffusion PDEs as well as extensive numerical experimentation for the linear elastostatics PDE. Since the test functions are sought in a
discontinuous polynomial space, using $r=p$ still result in a larger space than the trial as the 
discontinuous spaces contain additional degrees of freedom. Furthermore, in the limit $h \rightarrow 0$
the space $\VV$ is essentially $L^2$, i.e, any polynomial degree above constants is inherently 
an enrichment of the test space. 
\end{rem}

The approach to establishing AVS-FE approximations described until this point can be 
established in FE software with relative ease. However, there are other alternative interpretations 
of DPG methods which are even more straightforward in terms of implementation aspects.
The inventors of the DPG, Demkowicz and Gopalakrishnan refer to this as different "hats" of DPG methods~\cite{Demkowicz1} and the one we have explained here is that of a Petrov-Galerkin method 
with optimal test functions which leads to~\eqref{eq:FE_disc_weak}. As for the DPG method, we are also going to consider another "hat", in which a saddle point, or mixed, interpretation 
of the AVS-FE method which allows us to employ high level FE solvers such as FEniCS~\cite{alnaes2015fenics}.
Hence, let us introduce a new unknown function $\err$, the \emph{error representation} function. 
This function derives its name since it is a Riesz representer of the approximation error 
induced by the AVS-FE approximation~\eqref{eq:FE_disc_weak} of the weak formulation~\eqref{eq:weak_form}:
\begin{equation} \label{eq:err_rep}
\boxed{
\begin{array}{ll}
\text{Find }   \err \in \VV \; \text{ such that:}
\\[0.1in]
     \left(\,\err,(\vv,\www) \, \right)_\VV = F(\vv) -  B((\uu^h,\sig^h),(\vv,\www)), \\ &  \hspace*{-0.5in}\forall (\vv,\www)\in \VV.
 \end{array}}
\end{equation}
Again, the broken nature of the test space allows this function to be approximated on each element $\Kep$ \emph{a posteriori} to the solution of~\eqref{eq:FE_disc_weak} to be used as an error estimate 
and error indicator. 
Using basic arguments (see, e.g.,~\cite{Demkowicz1}) the following 
saddle point system can be established:
\begin{equation} \label{eq:saddle_point}
\boxed{
\begin{array}{ll}
\text{Find } (\uu^{\textcolor{black}{h}},\sig^{\textcolor{black}{h}}) \in \UUUh \wedge   \err \in \VV \; \text{ such that:}
\\[0.1in]
     \left(\,\err,(\vv,\www) \, \right)_\VV  -  B((\uu^h,\sig^h),(\vv,\www)) = - F(\vv), \\ &  \hspace*{-0.5in}\forall (\vv,\www)\in \VV, \\
     B((\yy^h,\zzz^h),\err ) \hspace*{1.35in} =  0, \\ &  \hspace*{-0.5in}\forall(\yy^h,\zzz^h)\in \UUUh.
 \end{array}}
\end{equation}

The discretization of $\err$ in~\eqref{eq:saddle_point} follows standard FE methodology and the 
space $\VV$ is to be discretized with discontinuous polynomials. Clearly, the \emph{global} computational cost of this saddle point system is larger than that of computing optimal test 
functions to establish~\eqref{eq:FE_disc_weak}. However, this cost is justified as the error 
representation function is to be used as an \emph{a posteriori} error estimator as well as element-wise
error indicators to be used in mesh adaptive strategies. Additionally, the effort in implementation 
into high level FE solvers with well established documentation and capabilities is embarrassingly 
low.
The norm equivalence~\eqref{eq:norm_equivalence} between the energy norm and the norm on $\VV$ 
of Riesz representers leads to the follwoing identity for the error representation function:
\begin{equation}
\label{eq:norm_equivalence_error}
\norm{(\uu- \uu^h,\sig-\sig^h)}{\text{B}} = \norm{\err}{\VV},
\end{equation}
which allows us to approximate the approximation error in the energy norm, which is not directly
computable due to the supremum, as well as element-wise error indicators:
\begin{equation} \label{eq:err_ind_est}
\boxed{
\begin{array}{ll}
\norm{(\uu- \uu^h,\sig-\sig^h)}{\text{B}} \approx \norm{\errH}{\VV},
\\[0.1in]
\eta = \norm{\errH}{\VVK},
 \end{array}}
\end{equation}
where $\errH$ is the approximation of $\err$ computed from the discretization of the saddle point 
system~\eqref{eq:saddle_point}.
\begin{rem}
We conclude this section by noting that the two interpretations of the AVS-FE in this section are 
completely equivalent. Hence, potential users that are limited by their available computer 
resources or software has the option of pursuing either interpretation being aware of their caveats. 
\end{rem}
%
%

%

\subsection{Error Estimates}
\label{sec:error_estimates}

In this section, we establish \emph{a priori} error estimates for the AVS-FE method. While we here 
assume that the components $\sig^h$ are discretized with continuous polynomials in $\SCZO$, the 
analysis can be performed with minor modifications using, e.g., RT or BDM discretizations. 
Furthermore, we assume that the optimal test functions, to be computed in the approach
of~\eqref{eq:FE_disc_weak} are sought in local polynomial spaces of the same degree as the trial 
functions. Equivalently, we assume that the discrete error representation function is in discontinuous 
polynomial spaces of the same degree as the trial functions. We shall use the arbitrary constant $C$
to denote generic mesh independent constants.

The starting point of our analysis is the best approximation property 
of the AVS-FE and DPG methods in terms of the energy norm~\cite{Demkowicz5}. Hence,
let $(\uu,\sig)\in \UUU$ be the exact solution of the weak formulation~\eqref{eq:weak_form}  
and $(\uu^h,\sig^h)\in \UUUh$ its approximation computed from the AVS-FE 
discretization~\eqref{eq:FE_disc_weak} or equivalently from the saddle point 
system~\eqref{eq:saddle_point}. Then, the
energy norm of the approximation error satisfies:
\begin{equation} \label{eq:best_approx}
\norm{(\uu-\uu^h, \sig - \sig^h)}{\text{B}} \leq \norm{(\uu-\vv^h, \sig - \www^h)}{\text{B}},
\end{equation}
where $(\vv^h,\www^h)$ are arbitrary functions in $\UUUh$. The proof of this inequality is established
using classical techniques from functional analysis and both inf-sup and continuity constants being
unity. Additionally, since the energy norm is an equivalent norm to $\norm{\cdot}{\UUU}$ 
on the trial space, we consequently
have the following quasi-best approximation property:
\begin{equation} \label{eq:Qbest_approx}
\norm{(\uu-\uu^h, \sig - \sig^h)}{\UUU} \leq C \,\norm{(\uu-\vv^h, \sig - \www^h)}{\UUU},
\end{equation}
where the mesh independent norm equivalence constant $C$ depends on the continuity constants of
of the bilinear form and a Fortin operator~\cite{nagaraj2017construction,demkowicz2020construction}.
Another key component in the following analysis is the convergence of polynomial interpolating functions. Hence,  there exist a global polynomial interpolation operator $\Pi_{hp}$   \cite{babuvska1987hp}:
\begin{equation} \label{eq:Pi_inter}
\ds  \Pi_{hp} \, : \, U \rightarrow  U^{hp}. 
\end{equation} 
Thus, $\Pi_{hp}(u)$ represents an interpolant of $u$ consisting of \textcolor{black}{globally continuous piecewise polynomials},  then~\cite{oden2012introduction}:
\begin{thm} \label{thm:classical_a_priori}
Let $u \in \SHR$ and $ \Pi_{hp}(u)  \in U^{hp}$ be the interpolant of $u$ \eqref{eq:Pi_inter}. Then, there exists $C > 0$
such that the interpolation error can be bounded as follows:
%
%
\begin{equation}
\ds  \norm{u-\Pi_{hp}(u)}{\SHS} \leq  C\, \frac{\ds h^{\,\mu-s}}{\ds p^{\,r-s}} \norm{u}{\SHR},
\end{equation} 
where $h$ is the maximum element diameter, $p$ the minimum polynomial degree of
interpolants in the mesh, $s \leq r$,  and $\mu = \rm{ min }$ $(p+1,r)$.
\end{thm}

To establish error estimates in terms of the energy norm, we first establish a bound on the Riesz 
representers of the trial functions, i.e., the optimal test functions. 
\begin{lem} \label{lem:energy_bound_riesz} 
Let $(\uu, \sig) \in \UUU$ be the exact solution of the AVS-FE weak formulation \eqref{eq:weak_form} and  $(\uu^h, \sig^h) \in \UUUh$ its corresponding AVS-FE approximation from~\eqref{eq:FE_disc_weak}. Then:
\begin{equation} \label{eq:energy_rate}
\norm{(\uu-\uu^h,\sig-\sig^h)}{\rm{B}} \leq   C\, \ds  \frac{ h^{\,\mu-1}}{ p_{\uu}^{\,r_{\pp}-1}},
\end{equation} 
where $h$ is the maximum element diameter, $\mu =$ $\rm{min}$ $(p_{\uu}+1,r_{\textcolor{black}{\pp}})$,  $p_{\uu}$ the minimum polynomial degree of approximation of $\uu^h$
in the mesh, and  $r_{\pp}$ the regularity of the solution $\pp$ of the distributional PDE underlying the Riesz representation problem~\eqref{eq:riesz_problem}
\end{lem} 
\emph{Proof:} The RHS of~\eqref{eq:best_approx} can be bounded by 
the error in the Riesz representers of the exact and approximate AVS-FE trial functions 
by the energy norm equivalence in~\eqref{eq:norm_equivalence}, and the map induced by the Riesz representation problem~\eqref{eq:riesz_problem} to yield:
\begin{equation} \label{eq:energy_bound}
\norm{(\uu-\uu^h,\sig-\sig^h)}{\rm{B}} \leq  \norm{(\pp-\pp^h,\rr-\rr^h)}{\VV},
\end{equation} 
where $(\pp, \rr) \in \VV$  are the exact Riesz representers of $(\uu, \sig)$ through~\eqref{eq:riesz_problem}, and $(\pp^h, \rr^h) \in \VVh$  are the approximate Riesz representers of $(\uu^h, \sig^h)$ through a FE discretization of~\eqref{eq:riesz_problem}.
The definition of $\norm{\cdot}{\VV}$ then gives:
\begin{equation*}
\begin{array}{lll}
\norm{(\uu-\uu^h,\sig-\sig^h)}{\rm{B}} \leq \textcolor{black}{ C ( 
\summa{\Kep}{} \{h_m \norm{\Nabla \pp -  \Nabla \pp^h}{\SLTO}\}
+  \norm{\pp-\pp^h}{\SLTO} } + \norm{\rr-\rr^h}{\SLTO} \textcolor{black}{)}. 
 \end{array}
\end{equation*}
\textcolor{black}{Since $\norm{\pp-\pp^h}{\SLTO} \leq \norm{\pp-\pp^h}{\SHOP} $ and $\norm{\Nabla \pp- \Nabla\pp^h}{\SLTO} \leq \norm{\pp-\pp^h}{\SHOP} $ , we get:
\begin{equation*}
\begin{array}{lll}
\norm{(\uu-\uu^h,\sig-\sig^h)}{\rm{B}} \leq \textcolor{black}{ C ( 
\summa{\Kep}{} \{h_m \norm{\pp-\pp^h}{\SHOP} }\}
+  \norm{\pp-\pp^h}{\SHOP}  + \norm{\rr-\rr^h}{\SLTO} \textcolor{black}{)}. 
 \end{array}
\end{equation*}
Now, we pick $h_m = h_{max} = h$, and 
}
trial functions that are polynomial interpolants for the Riesz representers $(\pp^h,\rr^h)$ of the same degree $p_{\uu}$. \textcolor{black}{Hence,} we bound $\norm{\rr-\rr^h}{\SLTO}$ using Theorem~\ref{thm:classical_a_priori} and $\norm{\pp-\pp^h}{\SHOP}$ by an extension of this theorem for broken Hilbert spaces as introduced 
by Rivi\'{e}re \emph{et al.}  in~\cite{riviere1999improved} to get:
\begin{equation}\label{eq:bound_temp}
\begin{array}{lll}
\norm{(\uu-\uu^h,\sig-\sig^h)}{\rm{B}} \leq \ds \textcolor{black}{C_1 h \frac{ h^{\,\mu_1-1}}{ p_{\uu}^{\,r_{\pp}-1}}}  \textcolor{black}{ + C_2} \ds \frac{ h^{\,\mu_1-1}}{ p_{\uu}^{\,r_{\pp}-1}} + \textcolor{black}{C_3} \ds \frac{ h^{\,\mu_2}}{ p_{\uu}^{\,r_{\rr}}}, 
 \end{array}
\end{equation}
where $\mu_1 =$ $\rm{min}$ $(p_{\uu}+1,r_{\pp})$, $\mu_2 =$ $\rm{min}$ $(p_{\uu}+1,r_{\rr})$, $r_{\pp}$, \textcolor{black}{$r_{\rr}$} the regularities of the Riesz representer\textcolor{black}{s} of the PDEs underlying~\eqref{eq:riesz_problem}. Since the \textcolor{black}{second} term in the RHS of~\eqref{eq:bound_temp} is  \textcolor{black}{dominant}, the proof is completed.
\\ \noindent $\square$

Next, with the the quasi-best approximation property~\eqref{eq:Qbest_approx}  at hand, we can readily
introduce  \emph{a priori} bounds in classical Sobolev norms. First, the bound in terms of the norm $\norm{\cdot}{\UUU}$ is governed by the following lemma:
\begin{lem} \label{lem:U_rates}
Let $(\uu, \sig) \in \UUU$ be the exact solution of the AVS-FE weak formulation~\eqref{eq:weak_form} and $(\uu^h, \sig^h) \in \UUUh$ its corresponding AVS-FE approximation through~\eqref{eq:FE_disc_weak}.
Then:
\begin{equation} \label{eq:U_rate}
\ds \exists \, C > 0 \, : \norm{(\uu-\uu^h,\sig-\sig^h)}{\UUU} \leq  C\, \left(  \,  \frac{ h^{\,\mu_1-1}}{ p_{\uu}^{\,r_{\uu}-1}} + \ds \,\frac{ h^{\,\mu_2-1}}{ p_{\uu}^{\,r_{\sig}-1}}  \right), 
\end{equation} 
where $h$ is the maximum element diameter, $\mu_1 =$ $\rm{min}$ $(p_{\uu}+1,r_{\uu})$,  $p_{\uu}$ the minimum polynomial degree of approximation of $\uu^h$, $\mu_2 =$ $\rm{min}$ $(p_{\uu}+1,r_{\sig})$,
in the mesh,  $r_{\uu}$ the regularity of the solution $\uu$ of the governing PDE~\eqref{eq:model_pde}, and $r_{\sig}$ the regularity of the solution $\sig$ of the governing first order system PDE~\eqref{eq:model_pde_first_order}. 
\end{lem} 
\emph{Proof:} By the quasi-best approximation property:
\begin{equation*}
\begin{array}{lll}
\ds \norm{(\uu-\uu^h,\sig-\sig^h)}{\UUU} & \leq C \, \norm{(\uu-\vv^h,\sig-\www^h)}{\UUU},
 \end{array}
\end{equation*}
the definition of the norm on $\UUU$~\eqref{eq:norms} leads to:
\begin{equation*}
\begin{array}{lll}
\ds \norm{(\uu-\uu^h,\sig-\sig^h)}{\UUU} & \leq C \, \{  \norm{\uu-\vv^h}{\SHOO}+\norm{\sig-\www^h}{\SHdivO} \}, 
 \end{array}
\end{equation*}
\textcolor{black}{since we use basis functions}
that are polynomial interpolants and note that $\norm{\sig-\www^h}{\SHdivO} \leq \norm{\sig-\www^h}{\SHOO}$, the approximation property in Theorem \ref{thm:classical_a_priori} with $s = 1$ gives:
\begin{equation*}
\begin{array}{lll}
\ds \norm{(\uu-\uu^h,\sig-\sig^h)}{\UUU} & \leq  \ds \left( C_1 \,  \frac{ h^{\,\mu_1-1}}{ p_{\uu}^{\,r_{\uu}-1}} + \ds C_2\,\frac{ h^{\,\mu_2-1}}{ p_{\sig}^{\,r_{\sig}-1}}  \right), 
 \end{array}
\end{equation*}
where we use the definitions of Lemma~\ref{lem:U_rates} for the $\mu$'s and $r$'s. Finally, combining the constants $C_1,C_2$ 
and noting that in the AVS-FE method we always pick $p_{\uu} = p_{\sig}$ the desired error bound is established. 
\\ \noindent $\square$

\begin{rem}
Note that choosing approximation spaces such as RT or BDM for the stress variable 
optimal error estimates can be established for these spaces. We refer to the text of Brezzi and Fortin~\cite{BrezziMixed} for details.
\end{rem}

\section{Numerical Verifications}
\label{sec:experiments}
To assess the performance of our method, we first consider a problem with a smooth solution which 
allows us to asses the convergence properties of the AVS-FE method to verify the \emph{a priori}
bounds of Section~\ref{sec:error_estimates}.
Then, we consider an
 example problem considered by Brenner in~\cite{brenner1993nonconforming},
with a manufactured exact solution that is dependent on the Poisson's ratio
which is used in a comparison between the AVS-FE method, the mixed FE method of Arnold and Winther~\cite{arnold2002mixed}, and the Bubnov-Galerkin FE method.
As  final numerical verifications, we present two engineering applications, $i)$ an example in which a 
commonly applied 
engineering structure, a cantilever beam, is considered and $ii)$ the deformation of a composite 
structure.     

In the numerical verifications presented in this section we use the FE solvers Firedrake~\cite{rathgeber2017firedrake} and FEniCS~\cite{alnaes2015fenics}. In particular, for all presented 
verifications using uniform meshes we use Firedrake, whereas in the case of mesh-adaptive refinements, 
we use FEniCS. In all cases, we use the linear solver MUMPS~\cite{amestoy2001fully} to perform the
inversion of the resulting stiffness matrices

\subsection{Asymptotic Convergence Studies}
\label{sec:convgence}

To present the convergence properties of the AVS-FE for nearly incompressible elastostatics, 
we consider a 2D model problem with a smooth exact solution which ensures the stress regularity 
is $H^1$, i.e., we use \textcolor{black}{$C^0$ continuous} approximations for both variables.
The domain is the unit square, i.e., $\Omega =  ( 0, 1 )\times ( 0,1 ) \subset \mathbb{R}^2$,
consisting of a material that is nearly incompressible with physical properties listed in Table~\ref{tab:material_data_incomp}, and a sinusoidal exact solution.
\begin{table}[t]
\centering
\caption{\label{tab:material_data_incomp} Material data for the nearly incompressible problem.}
\begin{tabular}{@{}lll@{}}
\toprule
{Property \hspace{6mm}} & {Symbol \hspace{6mm}} & {Value\hspace{6mm}}  \\
\midrule \midrule

Young's modulus & $E $ & $1500$ MPa     \\
Poisson's ratio & $\nu $  &  $0.4999$   \\

\bottomrule
\end{tabular}
\end{table}
\begin{equation}\label{eq:incompress_exact}  
\begin{array}{rl}
\ds \uu^{ex} (\xx) = \begin{Bmatrix}
           \ds u^{ex}_x (\xx) \\
           \ds u^{ex}_y (\xx) \\
\end{Bmatrix}  
= \begin{Bmatrix}
          \ds \text{sin} (\pi x) \, \text{sin} (\pi y)  \\
          \ds \text{sin} (\pi x) \, \text{sin} (\pi y)  \\
\end{Bmatrix}  .
 \end{array} 
\end{equation}
where $\lambda$ is the Lam\'e parameter:
Inspection of~\eqref{eq:incompress_exact} reveals the proper boundary conditions are homogeneous 
Dirichlet conditions on the entire boundary $\partial \Omega$ and the source $\ff$ is chosen such that 
it is the  differential operator  of the PDE~\eqref{eq:model_pde} acting on $\uu^{ex} (\xx)$.

\begin{figure}[h!]
\subfigure[ \label{fig:l2_err} $\norm{\uu-\uu^h}{\SLTO}$. ]{\centering
 \includegraphics[width=0.45\textwidth]{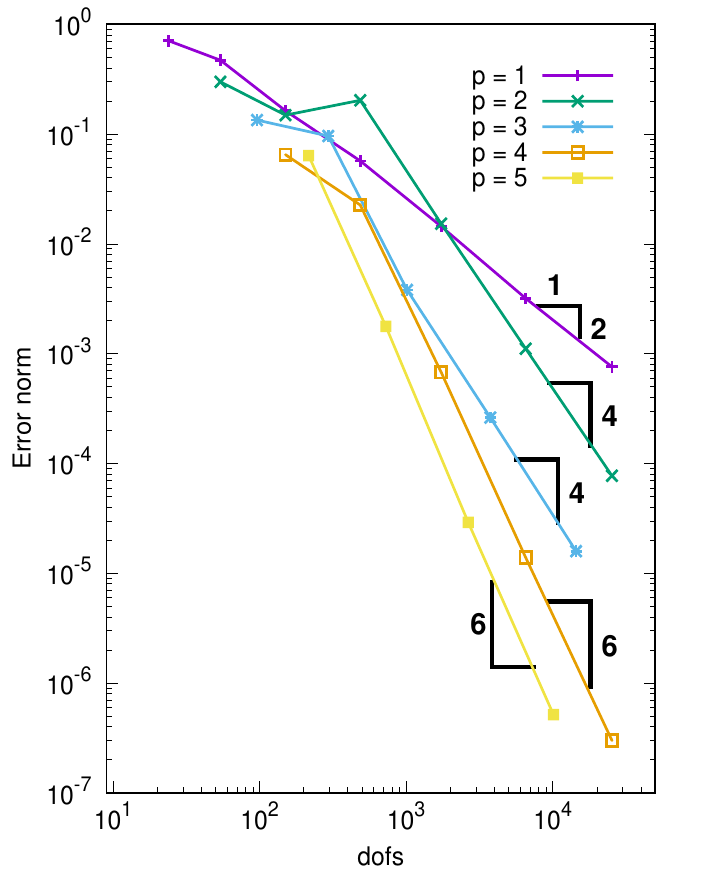}}
  \subfigure[ \label{fig:h1_err} $\norm{\uu-\uu^h}{\SHOO}$.]{\centering
 \includegraphics[width=0.45\textwidth]{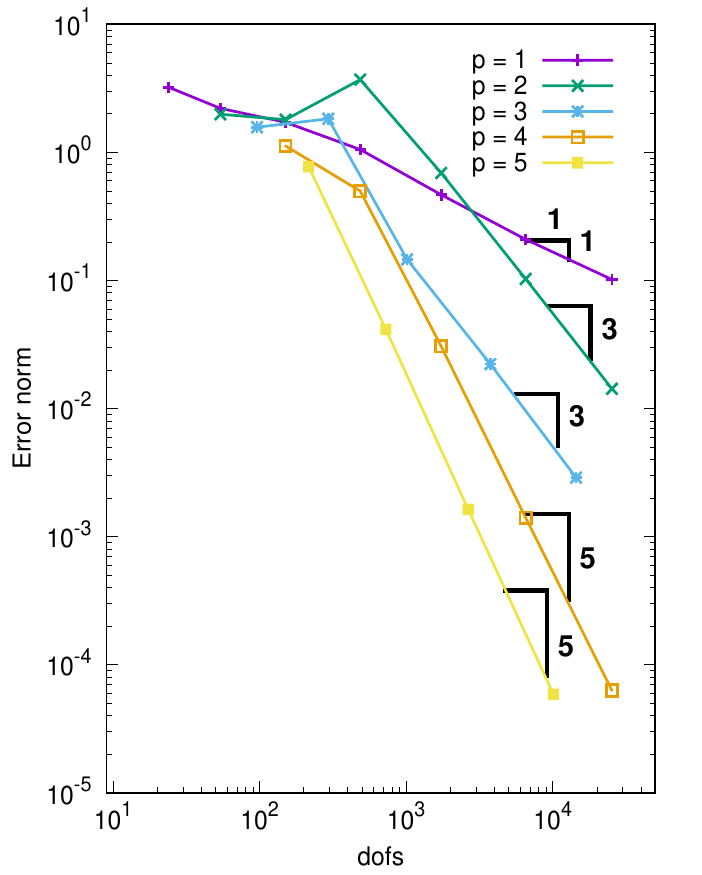}}
  \caption{\label{fig:l2h1} Asymptotic convergence results.}
\end{figure}

\begin{figure}[h!]
\subfigure[ \label{fig:UUU_err} $\norm{(\uu-\uu^h,\sig-\sig^h)}{\UUU}$. ]{\centering
 \includegraphics[width=0.45\textwidth]{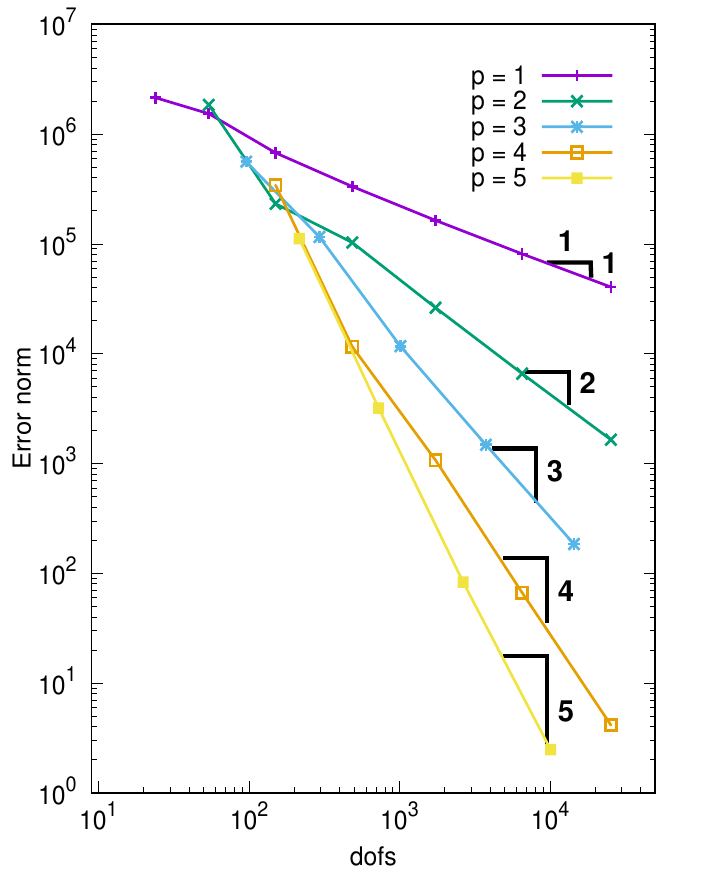}}
  \subfigure[ \label{fig:EE_Err} $\norm{(\uu-\uu^h,\sig-\sig^h)}{\text{B}}$.]{\centering
 \includegraphics[width=0.45\textwidth]{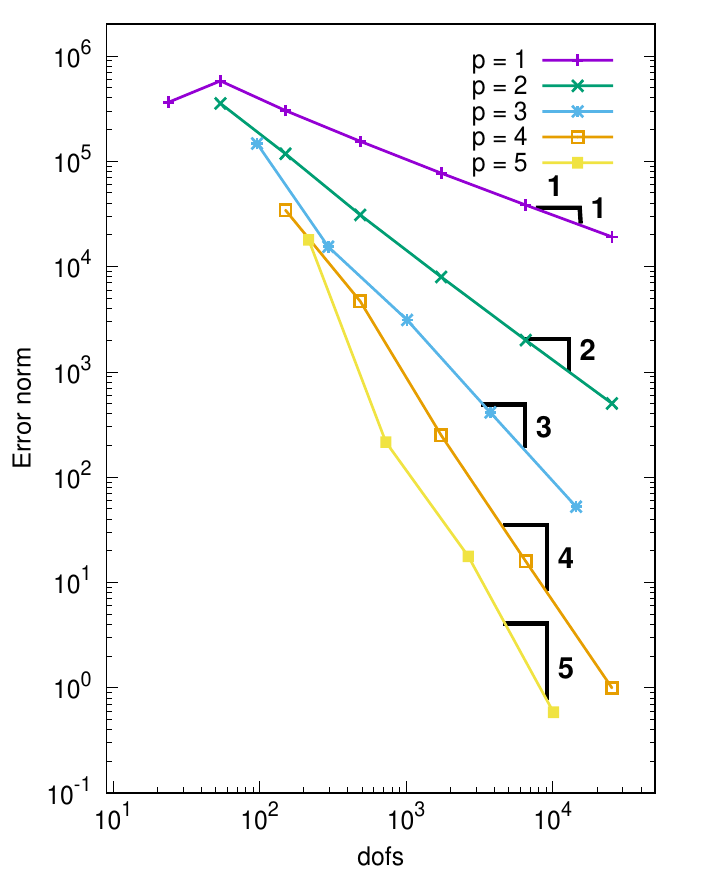}}
  \caption{\label{fig:UU_EE} Asymptotic convergence results.}
\end{figure}

As an initial verification, we consider uniform mesh refinements starting from a mesh consisting of two triangle elements,
and we use $\SCZO$ polynomials of increasing order. With this exact solution, the \emph{a priori}
error estimates from Section~\ref{sec:error_estimates} are:
\begin{equation} \label{eq:expected_rates}
\boxed{  
\begin{array}{rl}
 \ds \norm{(\uu-\uu^h,\sig-\sig^h)}{\UUU} \leq & C \, h^{p},
 \\[0.1in]
 \ds \norm{(\uu- \uu^h,\sig-\sig^h)}{\text{B}} \leq & C \, h^{p}.
 \end{array} }
\end{equation}
In Figures~\ref{fig:l2h1} and~\ref{fig:UU_EE}, we present the convergence history for the 
norms indicated in~\eqref{eq:expected_rates} with the exception of the $\norm{\sig-\sig^h}{\SHdivO}$ 
as this is a component of $\norm{(\uu-\uu^h,\qq-\qq^h)}{\UUU}$. In these figures we also show the 
errors in the individual norms $\norm{\uu-\uu^h}{\SLTO}$ and $\norm{\uu-\uu^h}{\SHOO}$.
The rates of convergence are as predicted in~\eqref{eq:expected_rates} for the energy norm and the norm on $\UUU$.
For the individual norms we observe the expected rates of convergence for polynomial FE approximations with the exception of the 
cases $p=4$ and $p = 2$ where the observed rates are one order higher in $\norm{\uu-\uu^h}{\SLTO}$ and
$\norm{\uu-\uu^h}{\SHOO}$.
This can be seen in Figures~\ref{fig:l2_err} and~\ref{fig:h1_err} where the slopes for $p=4$ and $p = 2$ are equal 
to the slopes for $p=5$ and $p = 3$, respectively.
This problem was also considered in~\cite{bramwell2012locking} for the DPG method and we point out that
the convergence rates of the AVS-FE approximations match the presented rates
in~\cite{bramwell2012locking}.

\subsection{Comparison with Other FE Methods}
\label{sec:conparisons}

To compare the AVS-FE method for nearly incompressible elasticity to other FE methods, 
we consider a model problem of which the solution is a function of the Poisson's ratio 
from~\cite{brenner1993nonconforming}.
The exact solution is given in~\eqref{eq:incompress_exact2}, we use the same Young's modulus as in the preceding example. However, we increase the Poisson's ratio 
to $0.49999999$ to make the problem more challenging with a Lam\'e parameter $\lambda$ of the order $10^{11}$. 
\begin{equation}\label{eq:incompress_exact2}  
\begin{array}{rl}
\ds \uu^{ex} (\xx) = \begin{Bmatrix}
           \ds u^{ex}_x (\xx) \\
           \ds u^{ex}_y (\xx) \\
\end{Bmatrix}  
= \begin{Bmatrix}
          \ds \text{sin} (2 \pi y) \left[ -1 + \text{cos} (2 \pi x) \right] +
          \frac{\ds \text{sin} (\pi x)\text{sin} ( \pi y)}{1+\lambda} \\
          \ds \text{sin} (2 \pi x) \left[ 1 - \text{cos} (2 \pi y) \right] +
          \frac{\ds \text{sin} (\pi x)\text{sin} ( \pi y)}{1+\lambda} \\
\end{Bmatrix}  ,
 \end{array} 
\end{equation}
The methods we consider are the Bubnov-Galerkin \textcolor{black}{(see Chapter 26 in~\cite{logg2012automated} for a description of the implementation of this method)} and \textcolor{black}{a mixed method of Arnold and Winther~\cite{carstensen2008arnold,arnold2002mixed}} in addition to
the AVS-FE method. For the Bubnov-Galerkin method, we approximate $\uu$ using $\SCZO$ linear
polynomials. Likewise, for the first-order system least squares method we use linear polynomials
 for both variables 
$\uu$ and $\sig$, and in the AVS-FE approximation we use linear polynomials and first order
 Raviart-Thomas bases for $\uu$ and $\sig$, respectively. \textcolor{black}{In the mixed method, $\uu$ is approximated with linear discontinuous polynomials and $\sig$ with the lowest order Arnold-Winther element (i.e., cubic).}
The initial mesh we consider is uniform, consisting of 8 triangular elements, and we proceed to perform uniform mesh refinements. 

\begin{figure}[h!]
\subfigure[ \label{fig:l2_comparison} $\norm{\uu-\uu^h}{\SLTO}$. ]{\centering
 \includegraphics[width=0.45\textwidth]{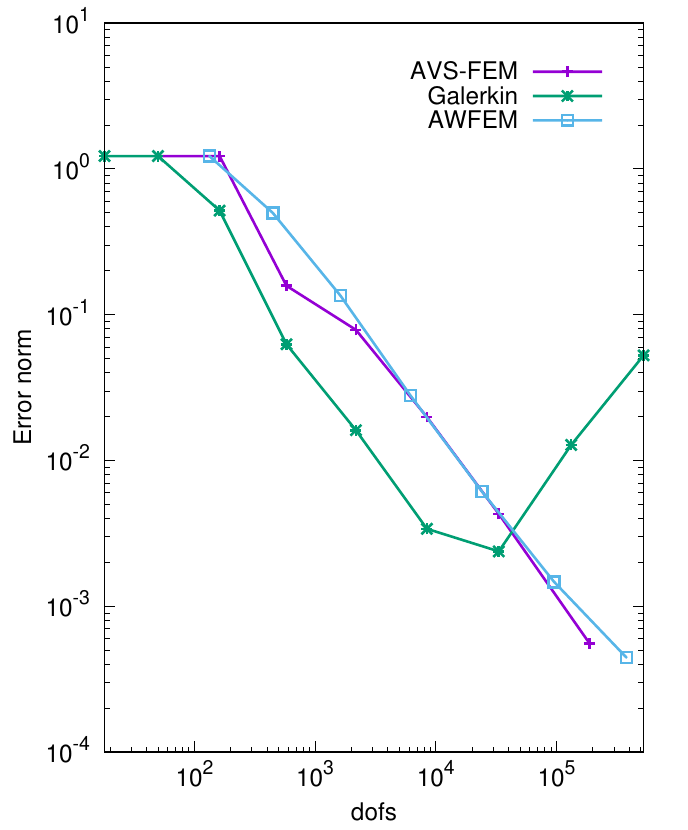}}
  \subfigure[ \label{fig:h1_comparison} $\norm{\uu-\uu^h}{\SHOO}$.]{\centering
 \includegraphics[width=0.45\textwidth]{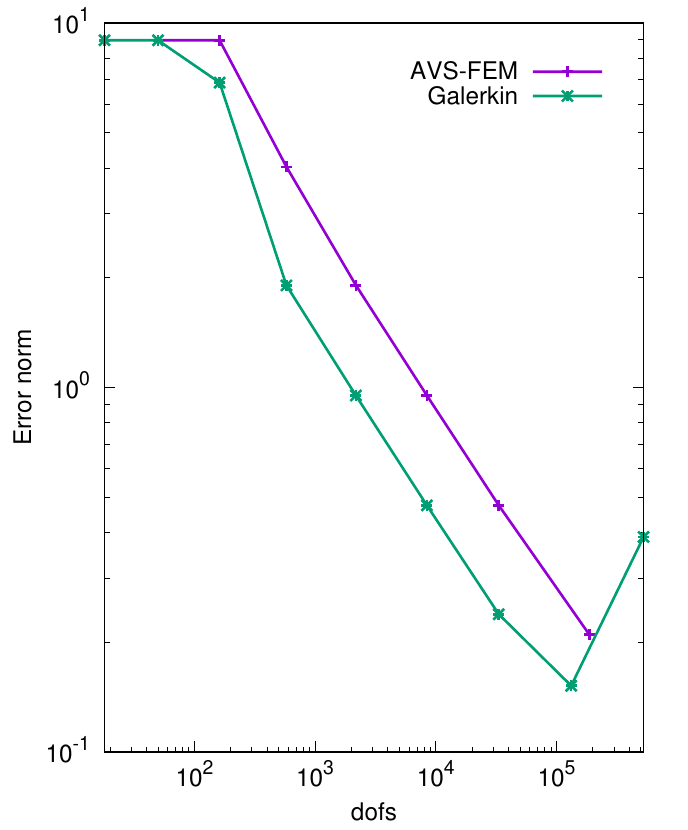}}
  \caption{\label{fig:convergence_results_incomp499999999} Comparison between the AVS-FE method, \textcolor{black}{the mixed method of Arnold and Winther,} and the Bubnov-Galerkin FE method.}
\end{figure}

In Figure~\ref{fig:convergence_results_incomp499999999} we compare the convergence history of these three methods in terms 
of the $L^2$, and $H^1$ norms \textcolor{black}{(only $L^2$ for the mixed method)} of the error on $\uu-\uu^h$. 
The Bubnov-Galerkin method performs well initially, but  upon continued  mesh refinements becomes unstable,
as indicated by the diverging errors. 
\textcolor{black}{The AVS-FE method
does not suffer from these issues, and retains the optimal rates of convergence in both
$L^2$ and $H^1$ norms}.
The increasing errors of the Bubnov-Galerkin FE method is likely 
due to ill conditioning of the resulting stiffness matrices. By lowering the Poisson
ratio this effect of ill conditioning is negated, as expected. Thereby limiting the applicability of
the Bubnov-Galerkin FE method for nearly incompressible materials.
\textcolor{black}{Finally, the mixed method of Arnold and Winther, remains stable as expected. However, in the last two refinements, the rate of convergence is reduced slightly. We attribute this to ill conditioning as in the Bubnov-Galerkin FE method. This has also been observed and studied by Carstensen \emph{et al.} in~\cite{carstensen2008arnold}, where it is noted that for pure Dirichlet boundary conditions the stiffness matrix condition number grows to infinity as $\nu \rightarrow 0.5$ (see Section 3.3 in~\cite{carstensen2008arnold}) ).   }

\subsection{Adaptive Mesh Refinement}
\label{sec:adaptivity}

Modern engineering materials often consist of multiple constituents, i.e., composites.  Material inclusions often lead to stress concentrations and need to be accurately 
resolved to ensure confidence in the engineering design. Resolution of such features is typically 
achieved by carefully constructed FE meshes which requires a significant effort from analysts, or 
mesh refinements.
Until this point, we have presented numerical verifications in which the mesh partitions are uniform
and their refinements are uniform as well, since the solutions are known to be rather well behaved, \emph{a priori}. 
While these certainly give us confidence in the AVS-FE approximations, the computational cost of 
uniform mesh refinements becomes very large as $h \rightarrow 0$. To this end, we consider adaptive 
mesh refinements that are guided by the built-in error indicators~\eqref{eq:err_ind_est}. 
As an adaptive strategy, we choose 
the marking strategy and refinement criteria of D{\"o}rfler~\cite{dorfler1996convergent} using 
a fixed parameter $\theta = 0.5$.

\begin{figure}[t]
\centering
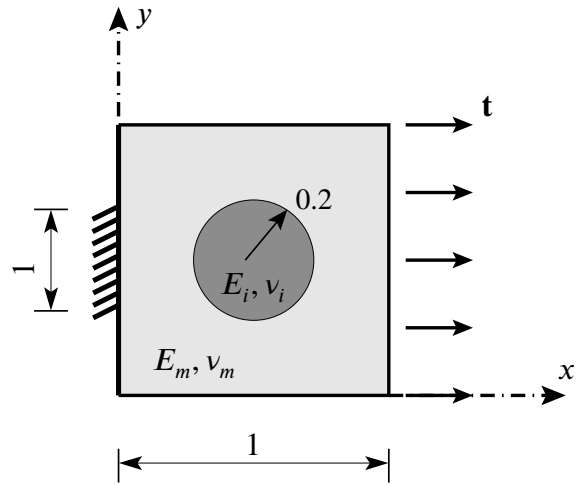
\caption{Linear elastic problem with inclusion.}
\label{fig:heterogeneous_solid}
\end{figure}

As an example of a simplified composite material we consider a unit square domain with an inclusion 
in the center at $x = y = 0.5$, as shown in Figure~\ref{fig:heterogeneous_solid}
The matrix material is nearly incompressible and the inclusion is a stiff material with Young's modulus 
several orders of magnitude larger than the bulk material. In Table~\ref{tab:mat_props} the properties
of the materials are listed. Materials with these properties are, e.g., Silicone based rubber for the 
matrix and an epoxy based inclusion.
\begin{table}[h]
\centering
\caption{\label{tab:mat_props}  Material properties.}
\begin{tabular}{@{}lll@{}}
\toprule
{Physical property} & {Young's Modulus \hspace{6mm}} & {Poisson's Ratio. \hspace{6mm}}  \\
\midrule \midrule

Matrix & 1,500 MPa & 0.49    \\
Inclusion & 10,000 MPa & 0.3   \\

\bottomrule
\end{tabular}
\end{table}
The boundary conditions are as shown in Figure~\ref{fig:heterogeneous_solid} i.e., a surface traction 
$\tr = \{100 \text{MPa},0 \}^T$ on the right side and a fully clamped boundary on a portion of the left side.
Hence, we expect significant concentration of stresses near the inclusion. We employ the error representation function~\eqref{eq:err_ind_est} as an \emph{a posteriori} error estimate. The goal of the adaptive algorithm 
is then to minimize this error based on its local error indicators.

\begin{figure}[h]
{\centering
\hspace{1in}\includegraphics[width=0.5\textwidth]{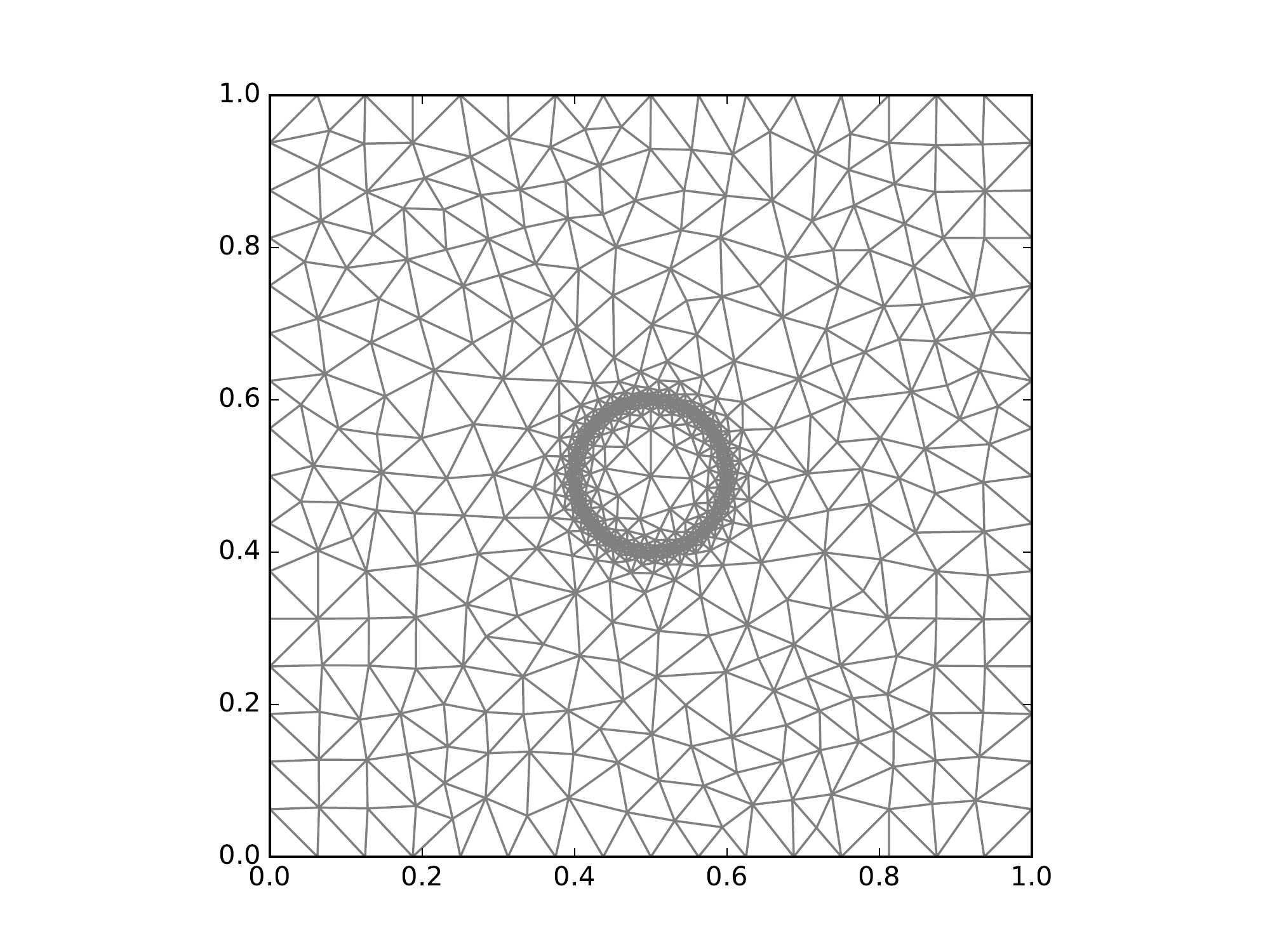}}
 \caption{\label{fig:in_mesh_inclusion} Initial mesh for the inclusion problem. }
\end{figure}

To establish an initial mesh taking into account the circular geometry with reasonable accuracy, 
we employ the built-in FEniCS~\cite{alnaes2015fenics} mesh generation functionality. In particular,
we use the tool "mshr" to create an inclusion consisting of 250 line segments to accurately represent
the circumference of the circle. The initial mesh is shown in Figure~\ref{fig:in_mesh_inclusion}.
For the AVS-FE approximations, we consider second order approximations for all variables in both cases
of RT and $\SCZO$ stress approximations. 

\begin{figure}[h!]
\subfigure[ \label{fig:sig_xx_adapt} $\sigma_{xx}^h(\xx)$. ]{\centering
 \includegraphics[width=0.45\textwidth]{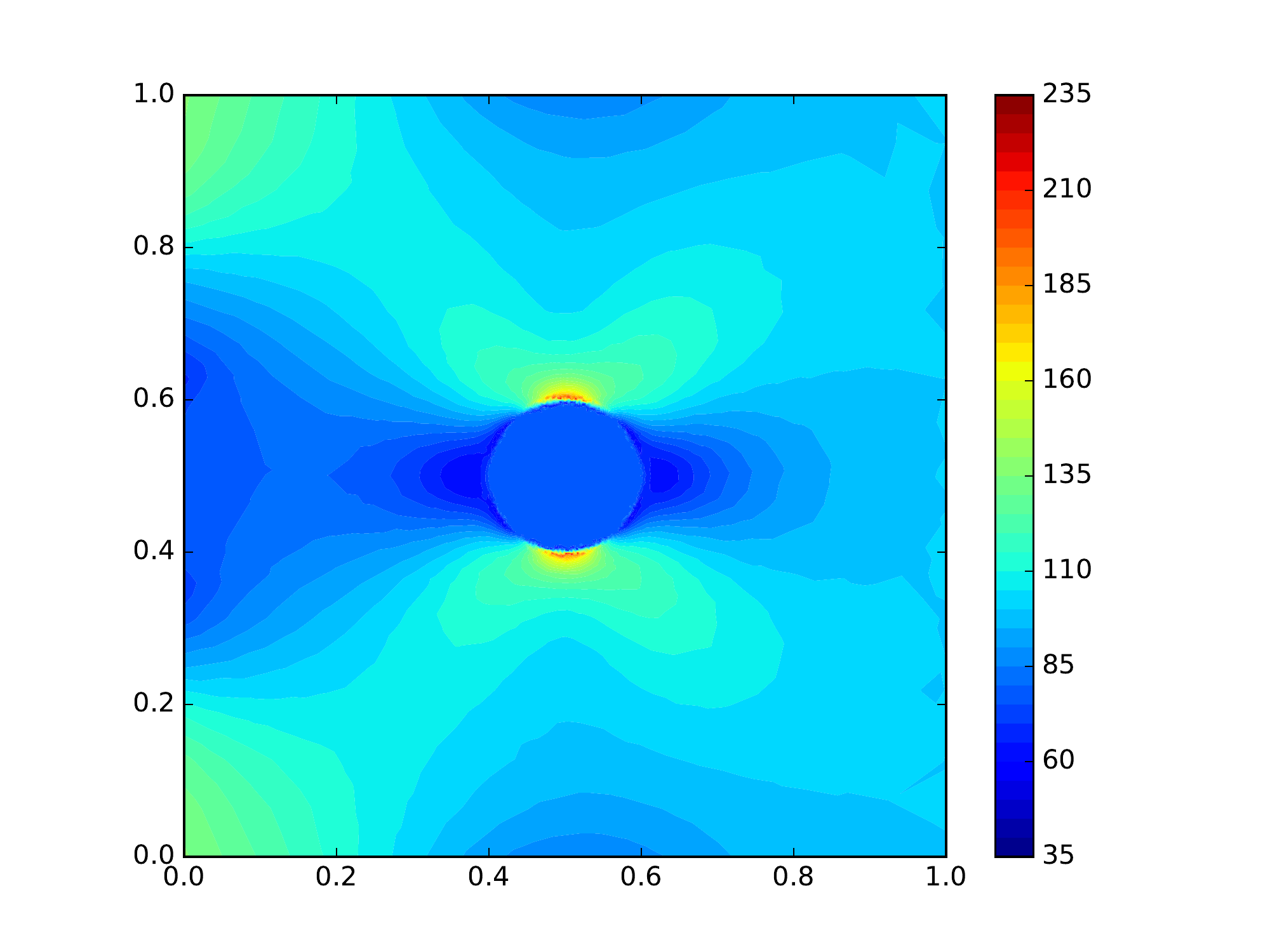}}
  \subfigure[ \label{fig:sig_yy_adapt} $\sigma_{yy}^h(\xx)$.]{\centering
 \includegraphics[width=0.45\textwidth]{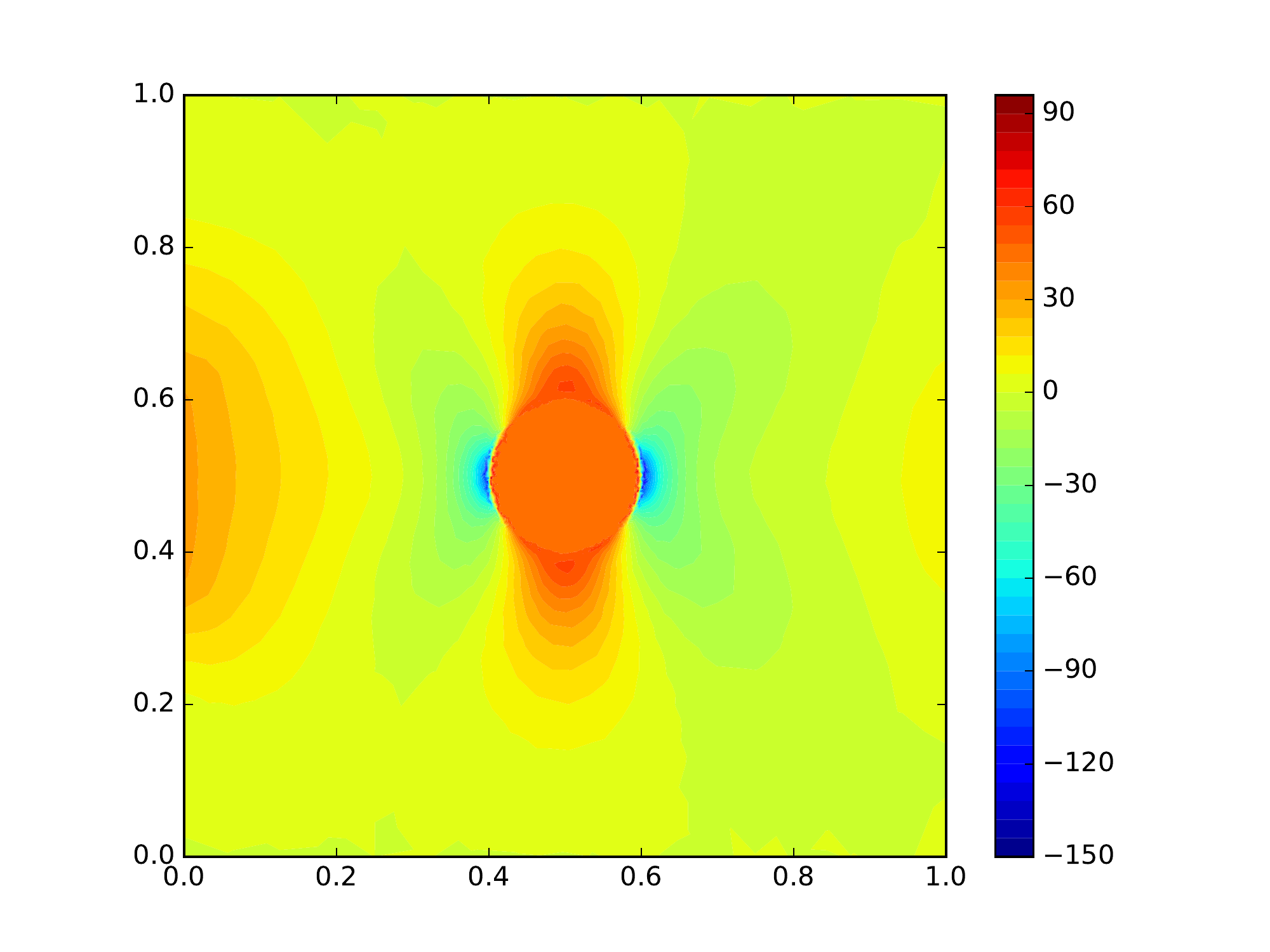}}
  \caption{\label{fig:final_stress_adaptive} AVS-FE approximate normal stress components ($\text{MPa}$).}
\end{figure}

\begin{figure}[h!]
\subfigure[ \label{fig:tau_xy_adapt} $\tau_{xx}^h(\xx)$. ]{\centering
 \includegraphics[width=0.45\textwidth]{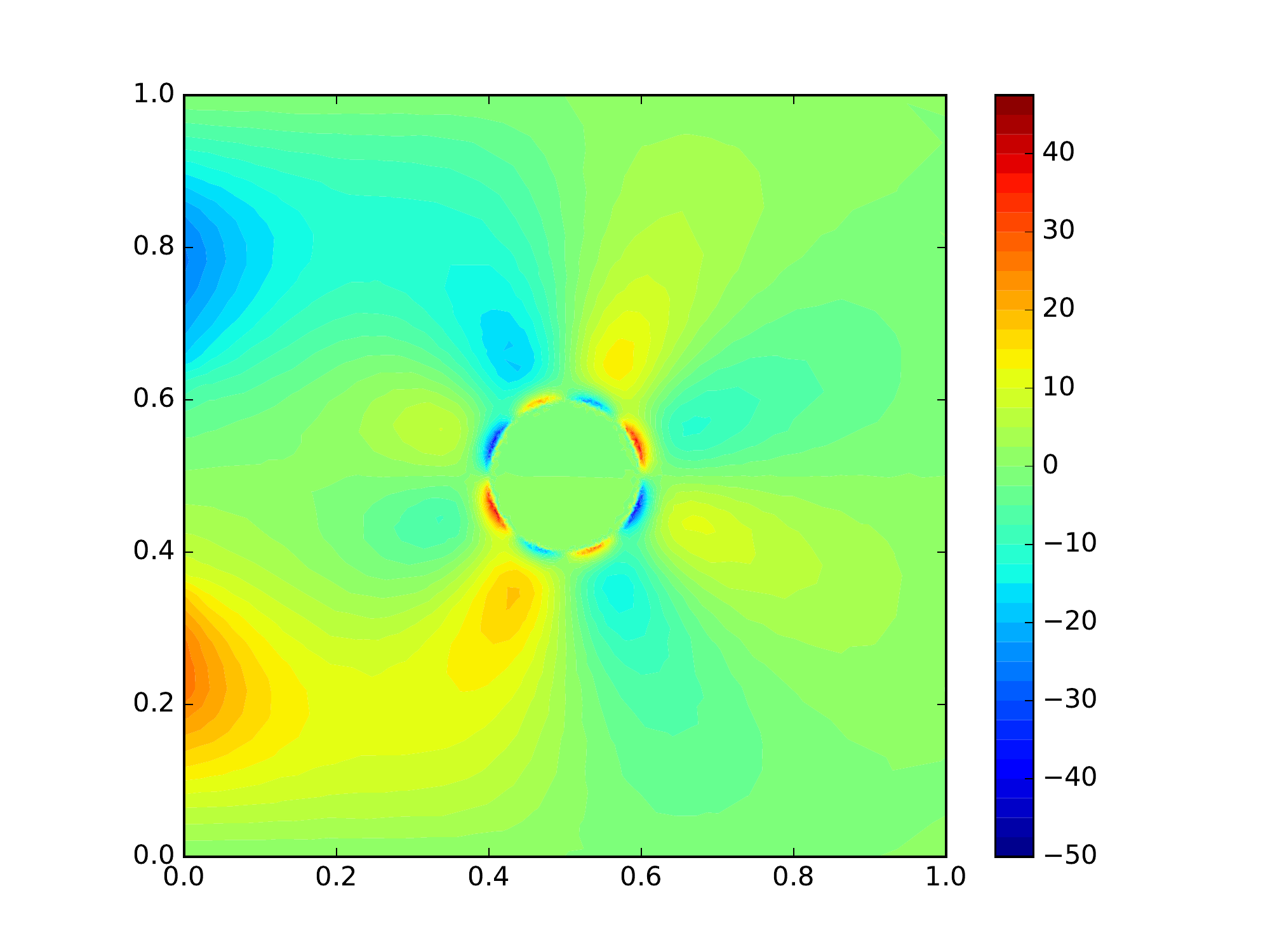}}
  \subfigure[ \label{fig:mesh_fin_Adap} Final adapted mesh using RT stress approximations.]{\centering
 \includegraphics[width=0.45\textwidth]{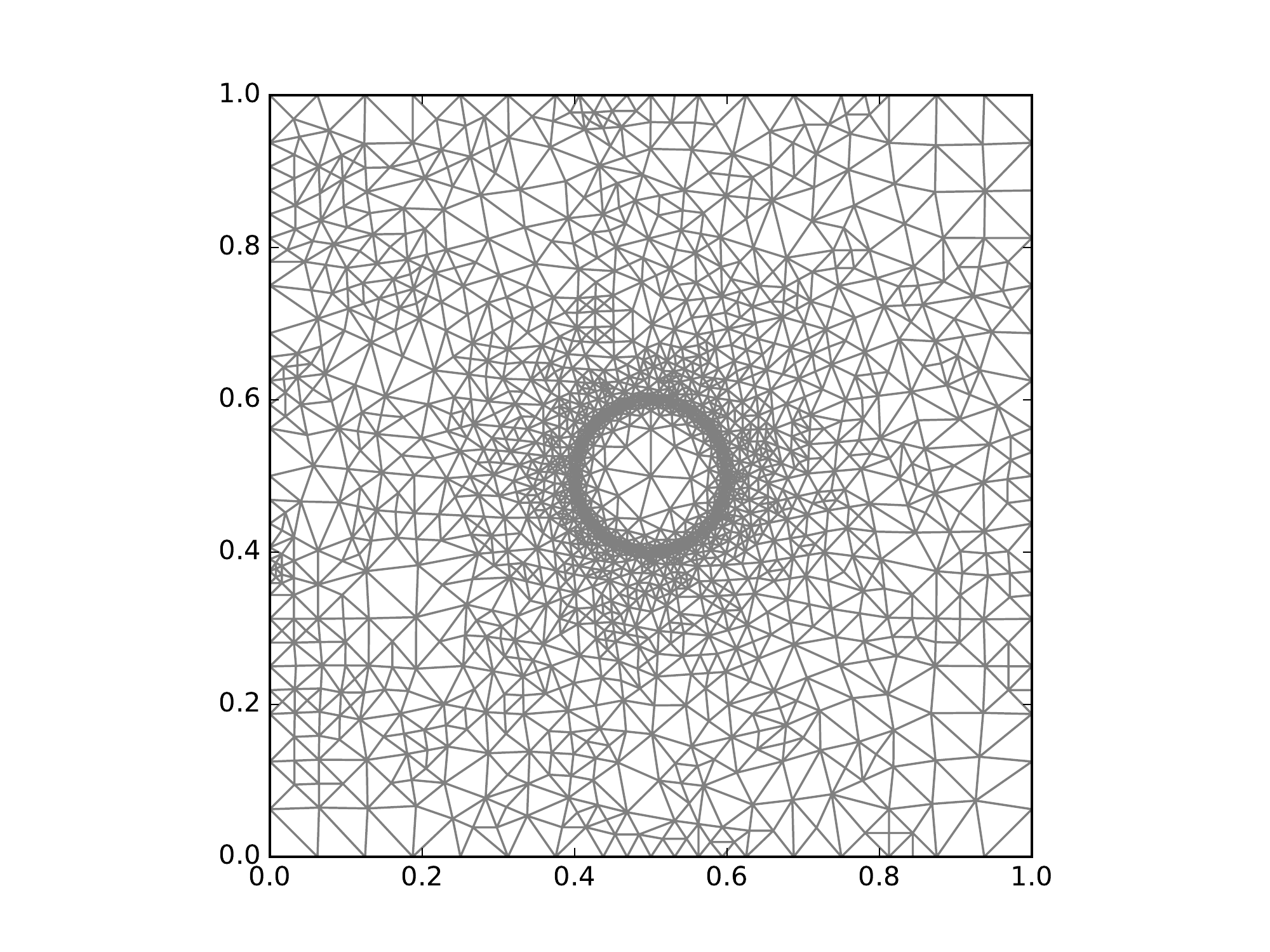}}
  \caption{\label{fig:final_shear_adaptive} AVS-FE approximate shear stress component ($\text{MPa}$) and final mesh.}
\end{figure}

The normal stresses $\sigma_{xx}$ and $\sigma_{yy}$ are presented in
 Figure~\ref{fig:final_stress_adaptive} where the stresses, 
as expected, are concentrated near the stiff inclusion. 
Furthermore, we see in this figure that the stresses perpendicular to the loading, i.e., $\sigma_{yy}$
are dominated by the Poisson effect which is another indication that the AVS-FE approximation 
of this problem is consistent with the expected physics. 
As a sanity check, we also note that the 
far field stresses are equal to the applied surface traction $\tr$, see Figure~\ref{fig:sig_xx_adapt}. 
The shear stress component $\tau_{xy}$ is shown in Figure~\ref{fig:final_shear_adaptive} along
 with the adapted mesh after 12 refinements.
The final adapted mesh shown in figure~\ref{fig:mesh_fin_Adap} highlights that the built-in error 
indicators as well as the marking strategy of D{\"o}rfler lead to mesh refinements in  
locations critical to proper resolution of physical features. 

Computing these results in the AVS-FE method is straightforward due to its built-in error 
estimate and discrete stability. Hence, the effort in implementation of the adaptive algorithm 
is minimal and requires only a FE solver with mesh refinement capabilities. 
However,
 it is worth mentioning that the required effort to establish similar 
results in the Bubnov-Galerkin FE method would be significant. While several \emph{a posteriori} error estimation techniques exists, the stability issue that arises with the jump in material coefficients 
will likely require significant efforts in analysis to ensure stability of the error estimation or 
prohibitively large mesh generation efforts.

\begin{figure}[h]
{\centering
\hspace{2in}\includegraphics[width=0.35\textwidth]{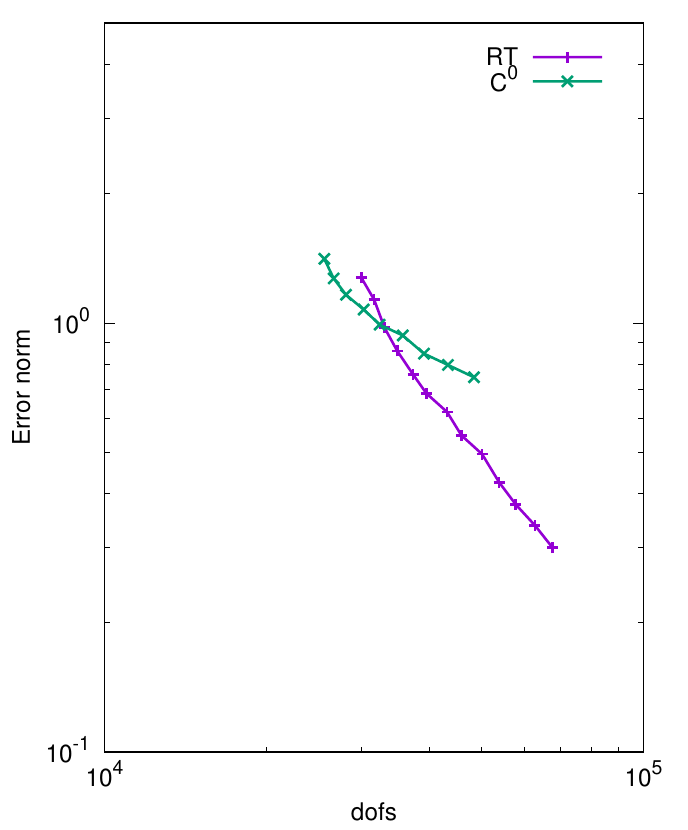}}
 \caption{\label{fig:c0_vs_rt} Convergence of the energy norm comparing $\SCZO$ and RT stress approximations \textcolor{black}{computed using the approximate energy norm, see~\eqref{eq:err_ind_est}}. }
\end{figure}

\begin{figure}[h]
{\centering
\hspace{1in}\includegraphics[width=0.5\textwidth]{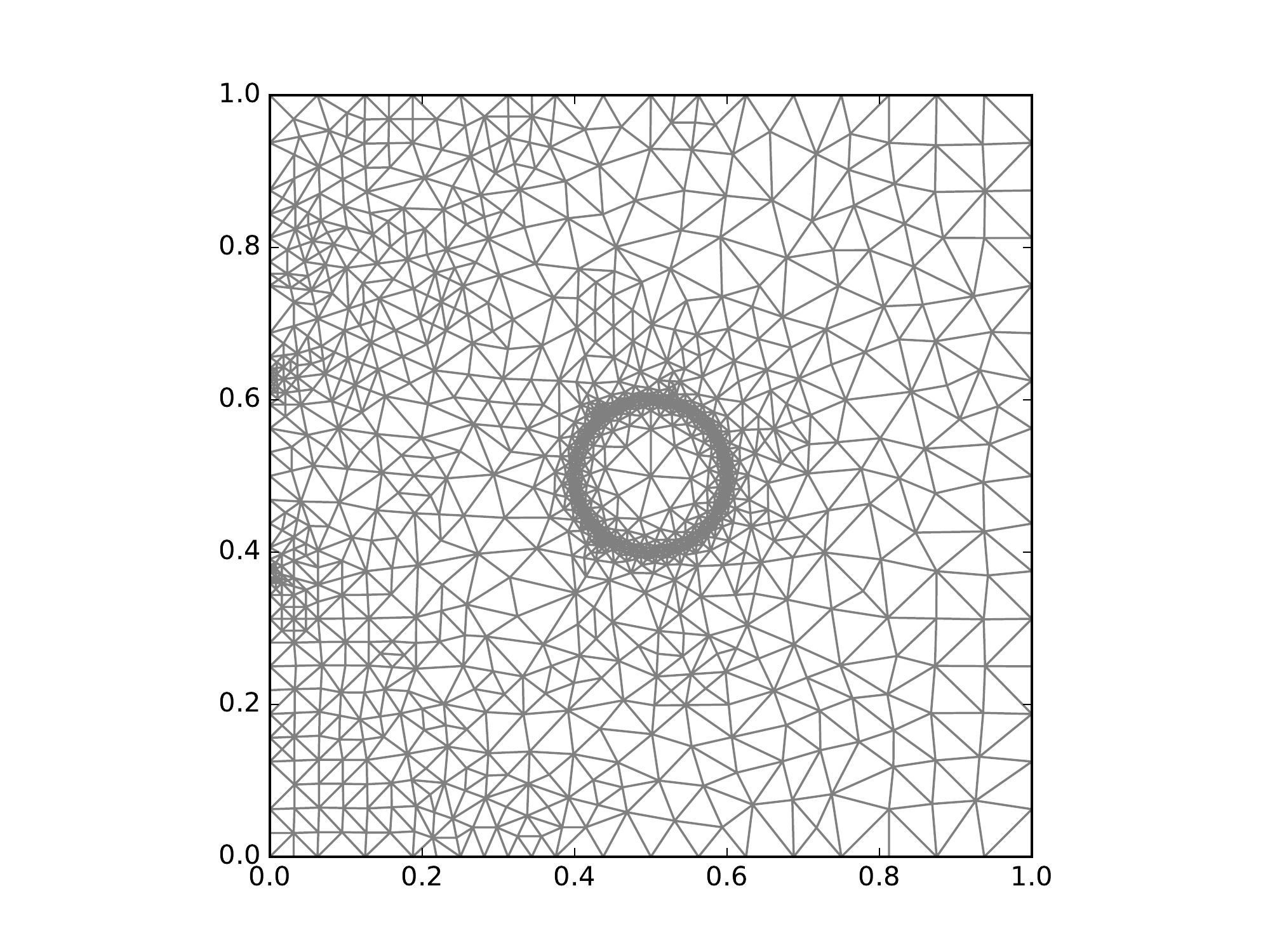}}
 \caption{\label{fig:fin_c0_mesh_inclusion} Final adapted mesh using $\SCZO$ stress approximations. }
\end{figure}

Finally, we compare $\SCZO$ and RT stress approximations. 
The inclusion leads to a stress field that has continuous normal components across the interface
between the two materials, whereas the tangential components are discontinuous. 
While this domain is still convex, we therefore expect that the generally advocated $\SCZO$ stress
approximations in the AVS-FE method 
will be less accurate than the standard mixed FE choice of RT discretization. Thus, as a 
verification we consider both types of stress approximations in this experiment and measure the 
difference between the two by the approximate energy norm. 
To this end, we employ the same adaptive algorithm for the case of $\SCZO$ stresses and 
 in Figure~\ref{fig:c0_vs_rt}, the convergence histories of the approximate error in the energy norm through\textcolor{black}{~\eqref{eq:err_ind_est}} for both cases are presented. 
As expected, in this case the RT approximations lead to lower energy norm errors due to the pollution 
effect of enforcing continuity of tangential stress components for the $\SCZO$ case. 
Finally, in Figure~\ref{fig:fin_c0_mesh_inclusion} the final adapted mesh for the case of $\SCZO$ stresses are presented. 
The overly restrictive stress approximations lead to mesh refinements that do not resolve the stress field near the inclusion in the same was as the RT case as evident from comparison of Figures~\ref{fig:mesh_fin_Adap} and~\ref{fig:fin_c0_mesh_inclusion}.
 Hence, the $\SCZO$ stress approximations are not applicable 
for this problem.

\subsection{Engineering Application: Non-Uniform Bending of a 2D Beam}
\label{sec:engineering_appl}
\begin{figure}[t]
\centering
\scalebox{0.75}{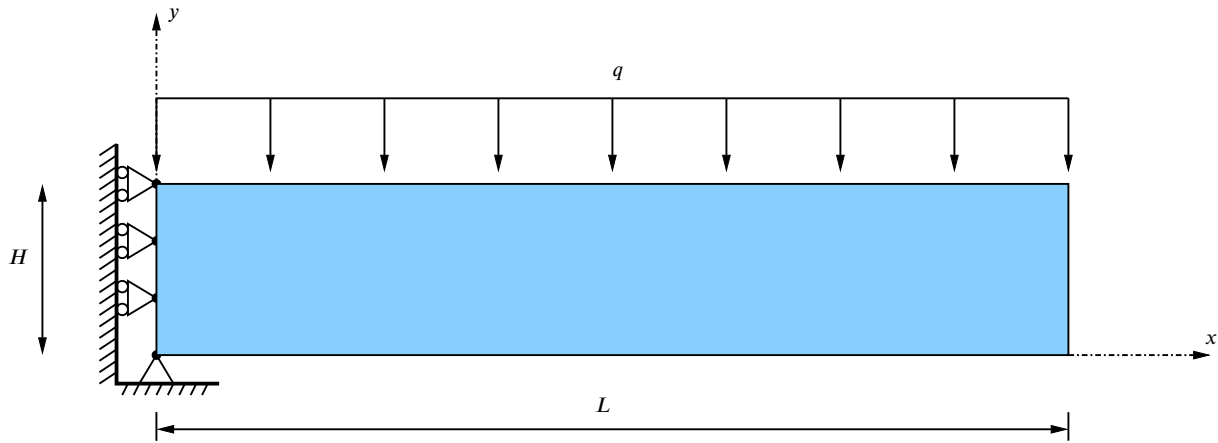}
\caption{Linear elastic beam problem.}
\label{fig:beam}
\end{figure}

To  further illustrate  the AVS-FE method in its application to FE analyses of nearly incompressible 
solids, we present a common engineering application of non-uniform bending  of beams.  The 2D problem 
concerns a beam with a length $L=2$m  and slenderness ratio $L/H=10$, and consists of a homogeneous 
linearly elastic isotropic material with a Young's Modulus equal to that of a nitrile based rubber, 
i.e., $E=1.5$MPa. The Poisson ratio  $\nu$ is chosen such that  $\frac{1}{2} - \nu = 10^{-9}$ and 
therefore the material is nearly incompressible.
The beam is subject to kinematic constraints along its left edge, where material particles are 
prohibited from moving in the $x$ direction but free to move in the $y$ direction. To eliminate the 
rigid body translation in the $y$ direction, the bottom left corner point is kept fixed.
In terms of loading, the beam is subject to a downward uniform distributed force $q$, of $3.33$N/m as 
depicted in Figure \ref{fig:beam}). 

\begin{figure}[h!]
 \includegraphics[width=.75\textwidth]{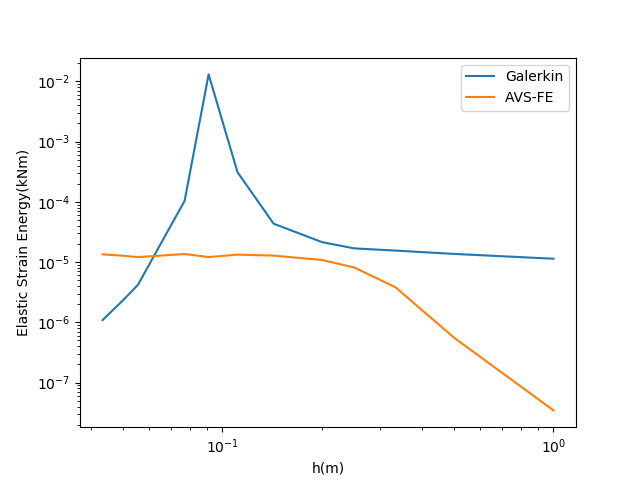}\centering
  \caption{\label{fig:ESE} Non-uniform bending of a nearly incompressible beam -- elastic strain energy evolution.}
\end{figure}

For our numerical experiment of the AVS-FE method, we start with an initial uniform mesh of four 
triangular elements (two elements in the length, versus one in the height) and subsequently conduct 
uniform $h$-refinements in which each element is partitioned into four new elements. In all 
computations, we employ quadratic (i.e., $p=2$) $C^0$ trial functions for the displacements and 
second order Raviart-Thomas functions 
functions for the stress variable. The evolution of the elastic strain energy of the 
AVS-FE solutions throughout the $h$-refinement process are shown in orange in Figure~\ref{fig:ESE}. In 
comparison, the corresponding results for the classical FE, or Bubnov-Galerkin, method are also shown in 
this graph in blue. These were established by using the same mesh partitions as for the AVS-FE method 
but employing the classical $C^0$ quadratic Lagrangian trial functions for the displacements and a 
standard displacement based weak formulation, as is 
common in classical Bubnov-Galerkin  analyses.

\begin{figure}[h!]
\subfigure[ \label{fig:NormalStress} Normal Stress $\sigma_{xx}$ (kPa). ]{\centering
 \includegraphics[trim={0cm 10cm 1cm 10cm},clip,width=1\textwidth]{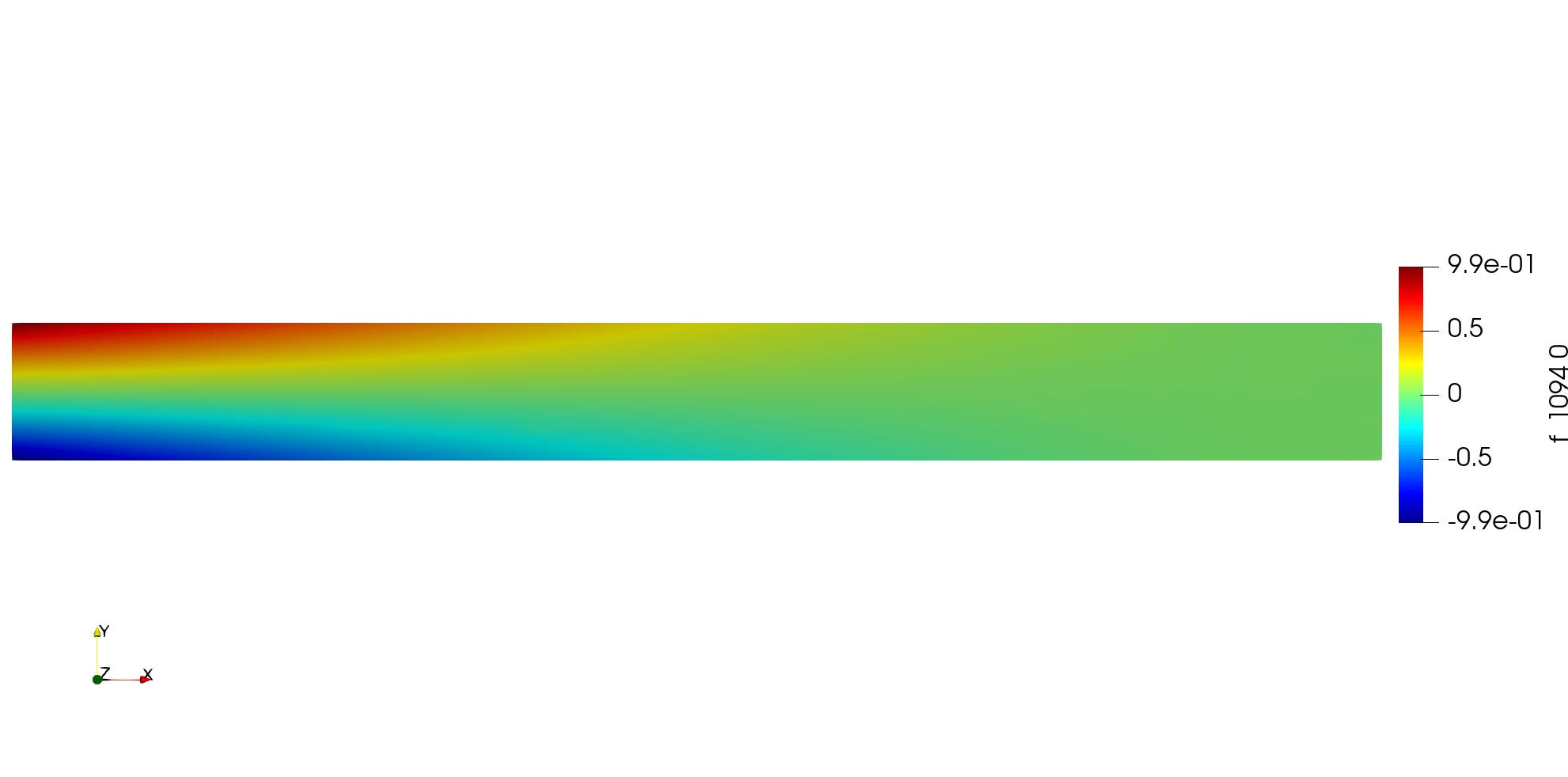}}
  \subfigure[ \label{fig:ShearStress} Shear stress $\tau_{xy}$ (kPa).]{\centering
 \includegraphics[trim={0cm 10cm 1cm 10cm},clip,width=1\textwidth]{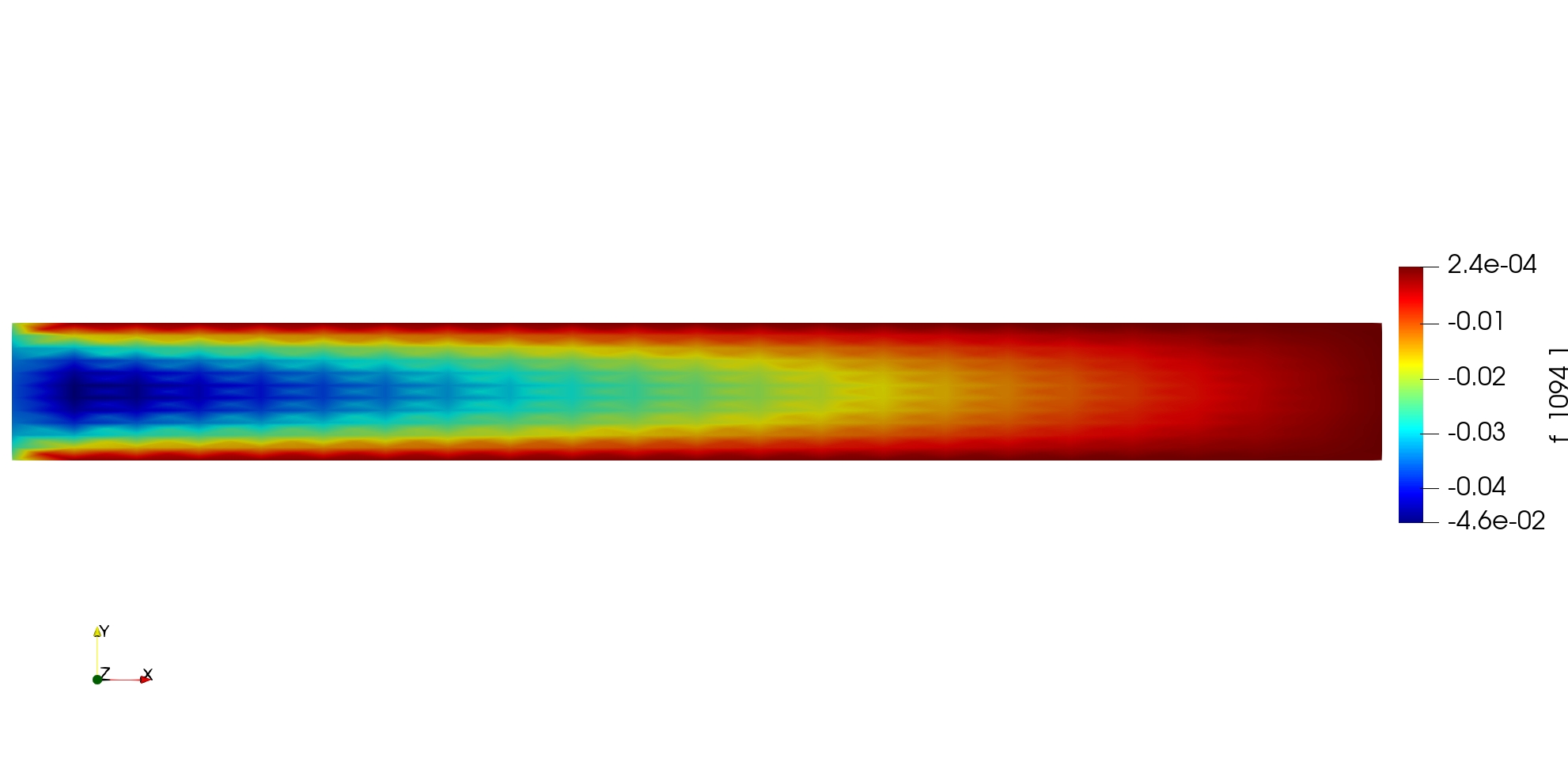}}
  \caption{\label{fig:Stress} AVS-FE stress distributions of beam in bending with $\frac{1}{2} - \nu = 10^{-9}$.}
\end{figure}

Figure~\ref{fig:ESE} shows that, as expected for this level of near incompressibility, the classical FE 
solutions fail to converge and start to exhibit spurious behavior as $h\longrightarrow 0$ with greatly 
changing values for the elastic strain energy between successive solutions. The AVS-FE method, on the 
other hand, maintains numerical stability from the onset and exhibits convergence upon continued uniform 
mesh refinements. The converged AVS-FE solutions are physically valid, as can be seen in 
Figure~\ref{fig:Stress}, in which contour plots of the distributions of the normal stress $\sigma_{xx}$ 
(Figure~\ref{fig:NormalStress}) and shear stress $\tau_{xy}$ (Figure~\ref{fig:ShearStress}) are shown of 
the AVS-FE solution for the mesh consisting of  $22\times 11$ elements. Since Raviart-Thomas trial functions have 
been used to compute the stress variables, minor discontinuities across inter-element edges can be 
observed, but these attenuate as the mesh is further refined. Hence, the AVS-FE method shows less 
sensitivity to the nearly incompressible constitutive behavior of the material than the classical FE 
method.

%
%
\section{Conclusions}
\label{sec:conclusions}

We have presented a mixed FE method that results in continuous FE approximations of both the displacement and (normal) stress fields in linear elasticity. 
In particular, we considered nearly incompressible materials as classical FE methods suffer 
from loss of discrete stability for such materials. 
The DPG philosophy of optimal test spaces leads to stable FE approximations without 
reformulation of the problem using the compliance tensor for the case for the nearly incompressible case.  
The DPG method we present distinguishes itself from the DPG method by only
breaking the test space while keeping the trial space globally confirming. 
Hence, in the corresponding FE approximation we use classical FE bases such as Lagrange polynomials and 
Raviart-Thomas functions. 

We present \emph{a priori} error bounds in terms of norms of the  numerical approximation
errors of both displacement and stress trial variables. These bounds are are all optimal in the sense 
that the approximation errors converge at rates equal to their underlying interpolating functions. 
Numerical verifications of the asymptotic convergence properties confirm the established 
error bound in all appropriate norms. A convergence study comparing our method to existing FE methods, 
i.e., Bubnov-Galerkin FE method and \textcolor{black}{the mixed FE method of Arnold and Winther~\cite{arnold2002mixed}}, show that the AVS-FE method is 
superior for the case of nearly incompressible materials for the presented verification. 
\textcolor{black}{In the verification presented, the Bubnov-Galerkin FE method suffered from a loss of stability as the mesh was refined, thereby resulting in loss of convergence. The mixed method of Arnold and Winther did not lose stability but for highly refined meshes a reduction in convergence rate was observed.} 
The AVS-FE method did not exhibit \textcolor{black}{these types of behavior}, nor have we observed such behavior for other verifications we have performed. 
We attribute this to the scaling term $h^2_m$ of the 
$H^1$ seminorm portion of $\norm{\cdot}{\VV}$~\eqref{eq:norms} as it ensures entries in the resulting 
stiffness matrix are of similar order of magnitude. 
However, we cannot rule out such behavior for the AVS-FE method for very fine meshes as this is an
issue of numerical linear algebra and not the AVS-FE method. 
Similar behavior has been observed by Storn in~\cite{storn2020relation} for the DPG method, 
where it is attributed to the inversion of the Gram matrix in the computation of optimal test
functions.

By considering a global saddle point form of the AVS-FE method, we establish both approximations of the 
displacement an stress fields as well as an approximation of an error representation function which 
measures the global energy error of the AVS-FE approximation. 
This error representation leads to \emph{a posteriori} error estimate\textcolor{black}{s} and error indicators 
which we employ in a mesh adaptive strategy. 
We present a numerical verification of a challenging physical application of a composite material 
where the built-in error indicator is used to drive adaptive mesh refinements to resolve 
the stress field in the composite.

While successful for the presented composite material, the built-in error indicator is a local 
indication of the residual of the AVS-FE approximation (see~\eqref{eq:norm_equivalence_error}). 
However, in certain applications, localized solution features may be of higher importance than the 
global energy error. Hence, goal-oriented error estimates and error indicators based on local 
quantities of interest~\cite{prudhomme1999goal} can provide alternative mesh refinement strategies. 
As shown in~\cite{valseth2020goal}, alternative AVS-FE goal-oriented error estimates are capable 
of accurately estimating these errors and driving goal-oriented mesh refinements. 
While we have considered a single AVS-FE weak formulation here~\eqref{eq:weak_form}, as 
mentioned in Remark~\ref{rem:alternative_forms}, this is not a unique choice. 
In~\cite{keith2016dpg}, multiple weak formulations for linear elasticity are considered for the 
DPG method, all of which provide slightly different FE approximations. Such investigation for the 
AVS-FE method and 
comparison of the DPG and AVS-FE methods for linear elasticity are postponed to future works.


\section*{Acknowledgements}
Authors Dawson and Valseth have been supported by the United States National Science Foundation - NSF PREEVENTS Track 2 Program, under  NSF Grant Number  1855047.
Authors  Kaul,  Romkes, and Valseth
have been supported by the United States National Science Foundation - NSF CBET Program,
under  NSF Grant Number 1805550. 

 \bibliography{references}

\end{document}